\documentclass[preprint,10pt]{elsarticle}
\usepackage{amssymb,amsmath,amsfonts,amsthm}
\usepackage{tikz,caption} 
\usetikzlibrary{positioning, arrows.meta, calc}
\usepackage{color}
\usepackage{latexsym}
\usepackage{geometry} 
\numberwithin{equation}{section}
\usepackage[normalem]{ulem} 

\newtheorem{theorem}{Theorem}[section]  
\newtheoremstyle{plainnoparen}
{3pt}{3pt}{\itshape}{}{\bfseries}{.}{.5em}
{\thmname{#1}\thmnumber{ #2}\thmnote{ {\normalfont #3}}}
\theoremstyle{plainnoparen}
\newtheorem{lemma}{Lemma}[section]
\theoremstyle{plain}
\newtheorem{corol}{Corollary}[section]
\newtheorem{prop}[theorem]{Proposition}

\newcommand{\p}{\partial}
\newcommand{\e}{\epsilon}
\newcommand{\n}{\nabla}
\newcommand{\dv}{{\rm div}}
\newcommand{\crl}{{\rm curl}}
\newcommand{\bu}{\mathbf{u}}
\newcommand{\bv}{\mathbf{v}}
\newcommand{\bn}{\mathbf{n}}
\newcommand{\mA}{\mathcal{A}}
\newcommand{\mC}{\mathcal{C}}

\newcommand{\mD}{\mathcal{D}}
\newcommand{\mV}{\mathcal{V}}
\newcommand{\mN}{\mathcal{N}}
\newcommand{\mP}{\mathbf{P}}
\newcommand{\mK}{\mathbf{K}}
\newcommand{\me}{\mathbf{e}}
\newcommand{\om}{\int_\Omega}

\newcommand{\2}{\|_{L^2}}
\newcommand{\3}{\|_{L^3}}
\newcommand{\4}{\|_{L^4}}

\newcommand{\0}{\|_{L^\infty}}
\newsavebox{\headpic}

\begin{document}
\begin{frontmatter}
	\journal{SIAM Journal on Mathematical Analysis}	
		
	\title{\textbf{The anelastic approximation and convergence rate for the compressible Navier--Stokes equations with ill-prepared data in bounded domains}}
		
	\author[1]{Yaobin Ou}
	\ead{ou@ruc.edu.cn}
	\author[1]{Ping Wang\corref{cor1}}
	\ead{wang-ping@ruc.edu.cn}
	\address[1]{School of Mathematics, Renmin University of China, Beijing 100872, P.R. China}
	\cortext[cor1]{Corresponding author.}
		
\begin{abstract}
    We justify the anelastic approximation and establish a convergence rate for local strong solutions to the isentropic compressible Navier–Stokes equations in a three-dimensional bounded domain, with ill-prepared initial data and vorticity-slip boundary conditions. When the Mach and Froude numbers vanish at the same rate $\e$, the system converges to an anelastic system with nonconstant hydrostatic density. Compared with the incompressible limit, the balance between singular pressure and gravity generates a variable-coefficient acoustic operator and additional commutators. To overcome these difficulties, we introduce the pressure perturbation adapted to the hydrostatic balance and establish weighted acoustic and vorticity estimates, yielding uniform-in-$\e$ estimates. We then use the weighted Helmholtz decomposition and time-averaging technique to prove the strong convergence of solenoidal component of the momentum and obtain an error bound controlled by $\e$ and the initial error in the projected momentum. Our analysis reveals the weighted cancellation mechanism that applies to singular limits around a class of nonconstant equilibria.
\vskip 8mm
\end{abstract}
	
		
\end{frontmatter}
\section{Introduction} 
	This paper is concerned with the anelastic approximation for viscous compressible flows in a gravitational field. We begin with the original compressible Navier--Stokes equations in a three-dimensional bounded domain $\Omega$ (see \cite{OP62}):
\begin{gather}\label{FOS_1}
	\rho_t+\dv(\rho \bu)=0,\\ \label{FOS_2}
	(\rho\bu)_t+\dv(\rho\bu\otimes\bu)+\n p=\dv \mathbb{S}(\bu)+\rho\n g,
\end{gather}
	where $\rho$, $\bu$, and $p=\rho^\gamma$ with $\gamma>1$ denote the density, velocity, and pressure, respectively. Moreover, $\mathbb{S}(\bu)=2\mu D(\bu)+\lambda\dv\bu I$ is the viscous stress tensor, with $D(\bu)=\frac{1}{2}(\n\bu+(\n\bu)^t)$, and $g=-x_3$ is the gravitational potential. The viscosity coefficients $\mu$ and $\lambda$ satisfy $\mu>0$ and $2\mu+3\lambda\geq0$.
	We non-dimensionalize the physical variables by scaling them by their reference values as follows (see Section 2.4 in \cite{Maj84})
\begin{gather*}
	\rho'=\frac{\rho}{r},\enspace p'=\frac{p}{\bar{p}},\enspace \bu'=\frac{\bu}{\bar{u}},\enspace g'=\frac{g}{\bar{g}},\\
	x'=\frac{x}{L},\enspace t'=\frac{t}{L/\bar{u}},\enspace \mu'=\frac{\mu}{\nu},\enspace \lambda'=\frac{\lambda}{\nu},
\end{gather*}
	where $r$, $\bar{p}$, $\bar{u}$, $L$, and $\nu$ are the reference density, pressure, velocity, length, and viscosity coefficient, respectively. Define the nondimensional constants
\begin{align*}
	\mathrm{Ma} = \frac{\bar{u}}{\sqrt{\gamma \bar{p}/r}}, \quad
	\mathrm{Fr} = \frac{\bar{u}}{\sqrt{\bar{g}L}}.
\end{align*}
    Dropping the primes and setting the remaining dimensionless  parameters to be 1, we obtain the nondimensional system (see \cite{FNS08, Mas07})
\begin{gather}\label{OS_1}
	\rho_t+\dv(\rho \bu)=0,\\ \label{OS_2}
	(\rho\bu)_t+\dv(\rho\bu\otimes\bu)+\frac{1}{\mathrm{Ma}^2}\n p=\dv \mathbb{S}(\bu)+\frac{1}{\mathrm{Fr}^2}\rho\n g.
\end{gather}
	
	We consider the strongly stratified regime ${\rm Ma}={\rm Fr}=\e$ and denote a solution to \eqref{OS_1}-\eqref{OS_2} by $(\rho^\e,\bu^\e)$. The hydrostatic profile is determined by the total mass $M=\om\rho_0^\e(x) dx$, which is assumed to be independent of $\e$. The background density $b$ is defined by the hydrostatic balance (see \cite{Mas07})
\begin{align}\label{hydro}
	\n(b^\gamma)-b\n g=0\enspace  {\rm and}\enspace \om b(x)dx=M \quad {\rm in}\enspace \Omega,
\end{align}
	which yields, for a constant $\Lambda$ uniquely determined below by the total mass $M$,
\begin{align}\label{defb}
	b=\left(\frac{\gamma-1}{\gamma}\right)^\frac{1}{\gamma-1}(\Lambda-x_3)_+^\frac{1}{\gamma-1}.
\end{align}
    In the sequel, we will assume that $b$ is bounded from below in $\Omega$. Indeed, since $\Omega$ is bounded, this no-vacuum condition holds whenever $\Lambda>\sup_{x\in\Omega}x_3$, which can be ensured by taking the total mass $M$ sufficiently large. Consequently, there exist constants $0<c_*<c^*<\infty$, depending only on $\Omega$, $\gamma$, and $M$, such that
\begin{align}\label{brho}
   	c_*\leq b(x)\leq c^*,\quad \|\n^k b^\alpha\|_{L^\infty(\Omega)}\leq C_{k,\alpha},
\end{align}
    for all nonnegative integers $k$ and all fixed $\alpha\in\mathbb{R}$. These bounds will be used crucially in the weighted estimates below, since the norms contain powers of $b^{-1}$ and other real powers of $b$.

	The hydrostatic balance leads to the formal limiting ansatz
\begin{align*}
	\rho^\e(t,x)\rightarrow b(x_3),\quad \bu^\e(t,x)\rightarrow\bv(t,x)\quad {\rm as}\enspace \e\to 0 .
\end{align*}
	Passing formally to the limit in the continuity equation \eqref{OS_1} yields the anelastic constraint 
\begin{align*}
	\dv(b\bv)=0.
\end{align*}
	Moreover, after subtracting the equilibrium relation \eqref{hydro}, the singular pressure-gravity contribution exhibits the weighted acoustic operator $b\nabla$ at leading order. In the limiting momentum equation, the constraint force takes the form $b\nabla\Pi$, where $\Pi$ is the Lagrange multiplier enforcing $\dv(bv)=0$ (see \cite{CotHol14}). The formal limiting system is therefore (see \cite{FNS08, Mas07}):
\begin{gather}\label{LS_1}
	\dv(b\bv)=0,\\ \label{LS_2}
	(b\bv)_t+\dv(b\bv\otimes\bv)+b\n\Pi=\dv \mathbb{S}(\bv).
\end{gather}
	
	The mathematical study of the low Mach number limit, with $1/{\rm Fr}$ fixed at $1$ or $0$, began with the work of Klainerman and Majda \cite{KlM81, KlM82}. Subsequent studies have treated the incompressible
	limit under a variety of boundary conditions and regularity assumptions. In the isentropic case, Lions and Masmoudi \cite{LMs98} used the group method to justify the low Mach limit for global weak solutions with ill-prepared data and various boundary conditions. By analyzing the damping of acoustic waves near the boundary, Desjardins et al. \cite{DGLM99} then proved strong convergence of the velocity field under Dirichlet boundary conditions. More recent works have considered strong solutions in bounded domains under different boundary conditions. Under general Navier-slip boundary conditions, Masmoudi, Rousset and Sun \cite{MRS22} used an anisotropic conormal framework to establish uniform regularity for the compressible Navier--Stokes system with ill-prepared data, while Ou and Yang \cite{OuY22} treated the case of vorticity-slip boundary conditions using a new weighted energy functional. Recently, Masmoudi et al. \cite{MSWZ25} investigated the Dirichlet problem and proved the incompressible and vanishing vertical viscosity limits by constructing approximate solutions that account for both fast waves and boundary layers. 
	
	In the non-isentropic case, the interaction between acoustic oscillations and the boundary is further complicated by temperature dynamics. Under Dirichlet boundary conditions, Jiang and Ou \cite{JOu11} obtained uniform estimates for well-prepared data without thermal conductivity. Jiang and Masmoudi \cite{JM15} refined the spectral boundary layer construction underlying the damping of acoustic modes and extended the result in \cite{DGLM99} to the full acoustic operator with heat conduction. With positive thermal conductivity and vorticity-slip boundary conditions, Ju and Ou \cite{JOu22} considered flows with large temperature variations, when the initial data are well-prepared in the sense of $L^2$-norm. Based on conormal Sobolev spaces, Sun \cite{Sun22} further established uniform regularity for Navier-slip boundary conditions and ill-prepared data. For more works concerning the low Mach number limit for compressible fluids, we refer the reader to \cite{Ala06, Fei07, Fei13}.
	
	In the regime where the Mach and Froude numbers vanish at the same rate, the pressure force and gravity remain balanced at the singular scale, and the limiting dynamics are governed by an anelastic system with a nonconstant hydrostatic density profile. This scaling is particularly relevant in modeling atmospheric flows (see, e.g., Ogura and Phillips \cite{OP62}), where high frequency acoustic waves are considered physically insignificant. However, a rigorous justification of the anelastic approximation is more delicate than the incompressible limit, since the large pressure gradient is coupled with the large gravitational force and the acoustic waves propagate around a stratified equilibrium. In a slab geometry with periodic boundary conditions, Feireisl et al. \cite{FNS08} investigated the anelastic approximation for global weak solutions to the compressible Navier--Stokes system. Meanwhile, Masmoudi \cite{Mas07} derived rigorously the anelastic approximation in general bounded domains by means of a local method. In exterior domains, Feireisl, Novotn\'y, and Petzeltov\'a \cite{FNP10} further verified the low Mach number limit for global weak solutions under strong stratification. For stratification generated by long-range potential forces, Feireisl \cite{Fei11} used spectral analysis and the RAGE theorem to prove the convergence to global weak solutions of the viscous anelastic system. Feireisl et al. \cite{FKNZ16} later studied global finite-energy weak solutions to the compressible Navier--Stokes system coupled to a renormalized transport equation for the potential temperature. They obtained convergence to an anelastic system coupled to a transport equation, for well-prepared data in the high Reynolds number regime on the flat torus and for ill-prepared data on an infinite slab. Combining the relative energy inequality with frequency localized Strichartz estimates, Donatelli and Feireisl \cite{DFi17} justified the anelastic approximation for global weak solutions in the high Reynolds number regime.
	
	More recently, the theory of singular limits under strong stratification has been extended to the cases of rotation, heat conduction, degenerate viscosity, and vacuum boundaries. Using a relative entropy method, Bocchi, Fanelli and Prange \cite{BFP22} connected fast rotating barotropic flows with anisotropic viscosity and stratification to quasi-geostrophic dynamics with Ekman boundary layers. Later on, Basari\'c et al. \cite{BBFOT23} studied a strongly stratified, heat-conducting fluid confined between parallel plates and identified Majda's layered ``stack-of-pancake'' flow as the limit dynamics. For degenerate, density-dependent viscosities, Fanelli and Zatorska \cite{FaZ23} proved convergence to a generalized anelastic system, while Chaudhuri et al. \cite{CFLZ25} removed the cold-pressure term for standard isentropic pressure laws. For the one-dimensional vacuum free boundary problem, Ou \cite{Ou21} established uniform estimates for global classical solutions and showed that both the combined singular limit and the large time asymptotic behavior converge to the same steady state. These results largely concern weak solutions and exploit specific geometries or additional physical structures. For strong solutions in multi-dimensional bounded domains, the anelastic approximation also requires higher-order estimates that remain uniform in $\e$.

	In this paper, we rigorously justify the anelastic approximation for isentropic compressible Navier--Stokes equations with ill-prepared initial data and slip boundary conditions in three-dimensional bounded domains. In the finite-energy weak solution frameworks of Feireisl et al. \cite{FNS08} and Masmoudi \cite{Mas07}, a central difficulty is passing to the limit in the nonlinear terms in the presence of fast acoustic oscillations. For strong solutions, uniform higher-order estimates are also needed. The nonconstant hydrostatic density produces singular commutators under spatial differentiation of the acoustic equations. We obtain uniform $H^2$-bounds and weighted higher-order estimates, and then prove strong convergence of the projected momentum to that of the limiting anelastic flow.

\subsection{Main results}		
    The acoustic dynamics are governed by the nonconstant hydrostatic profile $b$. We therefore introduce a pressure perturbation adapted to the pressure–gravity balance, which exposes the weighted acoustic structure.

    We first introduce the \textit{density perturbation} $\phi^\e$ and the \textit{pressure remainder} $\pi^\e$ through the expansion
\begin{align}\label{rhop}
	\rho^\e=b+\e\phi^\e,\quad p^\e=(\rho^\e)^\gamma=b^\gamma+\gamma\e b^{\gamma-1}\phi^\e+\e^2\pi^\e,
\end{align}
    where the \textit{pressure remainder} $\pi^\e$ is defined exactly by
\begin{align}
	\pi^\e:=\frac{(\rho^\e)^\gamma-b^\gamma-\gamma\e b^{\gamma-1}\phi^\e}{\e^2}.
\end{align}
    With the hydrostatic balance \eqref{hydro}, a direct computation (see also \cite{Mas07}) yields
\begin{align*}
	\frac{1}{\e^2}(\n p^\e-\rho^\e\n g)=&\frac{1}{\e^2}\left(\n b^\gamma+\gamma\e\n(b^{\gamma-1}\phi^\e)-(b+\e\phi^\e)\n g\right)+\n\pi^\e\\
	=&\frac{1}{\e}\left(\gamma\phi^\e\n(b^{\gamma-1})+\gamma b^{\gamma-1}\n\phi^\e-\phi^\e\n g\right)+\n\pi^\e\\
	=&\frac{1}{\e}\left(\gamma\phi^\e\frac{\gamma-2}{\gamma-1}\n(b^{\gamma-1})+\gamma b^{\gamma-1}\n\phi^\e\right)+\n\pi^\e\\
	=&\frac{\gamma}{\e}b\n(b^{\gamma-2}\phi^\e)+\n\pi^\e.
\end{align*}
    By introducing the \textit{modified pressure perturbation} 
\begin{align}\label{defq}
	q^\e:=\gamma b^{\gamma-2}\phi^\e,
\end{align}
    the pressure and density can be rewritten as
\begin{align}\label{defp}
	p^\e=b^\gamma+\e bq^\e+\e^2\pi^\e,\quad \rho^\e=(b^\gamma+\e bq^\e+\e^2\pi^\e)^\frac{1}{\gamma}.
\end{align}
    Setting
\begin{align}
	\mA(b)=\frac{1}{\gamma}b^{2-\gamma},\quad \mD_b\bu^\e:=b^{-1}\dv(b\bu^\e),
\end{align}
    we obtain from \eqref{defq}
\begin{align}\label{defr}
	\rho^\e=b+\e\mA(b)q^\e.
\end{align}
    Therefore, we can rewrite the nondimensional system \eqref{OS_1}-\eqref{OS_2} as follows:
\begin{gather}\label{ts_1}
	\mA(b)q^\e_t+\dv(\mA(b)q^\e\bu^\e)+\frac{1}{\e}\dv(b\bu^\e)=0 \quad\quad {\rm in}\enspace (0,T]\times\Omega,\\ \label{ts_2}
	\rho^\e(\bu_t^\e+\bu^\e\cdot\n \bu^\e)+\frac{1}{\e}b\n q^\e+\n\pi^\e=\mu\Delta \bu^\e+(\mu+\lambda)\n\dv \bu^\e \quad\quad {\rm in}\enspace (0,T]\times\Omega.
\end{gather}
    We impose the following initial and boundary conditions:
\begin{gather}\label{TS_ic}
	(q^\e,\bu^\e)(0,x)=(q_0^\e,\bu_0^\e)(x) \quad\quad {\rm in}\enspace \Omega,\\ \label{TS_bc}
	\bu^\e\cdot \bn=0,\quad \bn\times\crl \bu^\e=0 \quad\quad {\rm on}\enspace (0,T]\times\p\Omega,
\end{gather}
    where $\bn$ is the unit outward normal vector to the smooth boundary $\p\Omega$.
    
    \noindent\textbf{Remark 1.1.} Note that the second condition in \eqref{TS_bc}, namely the vorticity-slip boundary condition, is closely related to the Navier-slip boundary condition (see \cite{Zaj98})
\begin{align*}
	\tau\cdot(2\mathbb{S}(\bu^\e)\bn)+\alpha(x)\bu^\e\cdot\tau=0\quad {\rm on}\enspace (0,T]\times\p\Omega,
\end{align*} 
    where $\tau$ denotes the unit tangential vector along $\p\Omega$, and $\alpha(x)$ is a given function. The link between these two boundary conditions relies on the identity (see \cite{XX13})
\begin{align*}
	\tau\cdot(2\mathbb{S}(\bu^\e)\bn-\bn\times\crl\bu^\e)=-2\tau\cdot(\mathbb{S}(\bn)\bu^\e)\quad {\rm on}\enspace \p\Omega,
\end{align*}
    provided that $\bu^\e\in H^2(\Omega)$ with $\bu^\e\cdot \bn=0$ on $\p\Omega$.
    
    The anelastic approximation requires estimates for local strong solutions that are uniform in $\e$. The following energy functional controls the acoustic variables in weighted norms. Since $\dv\bu^\e=\mD_b\bu^\e-\bu^\e\cdot\n\log b$, estimates for $\mD_b\bu^\e$ yield the corresponding estimates for $\dv\bu^\e$, once the lower-order velocity terms are controlled. 
    
    \noindent\textbf{Definition 1.1.} We define the weighted energy functional
\begin{align*}
	M^\e(t):=\sup_{s\in[0,t]}E^\e(s)+\int_0^tD^\e(s)ds,
\end{align*}
    where
\begin{align*}
	E^\e(t):=&\|(q^\e,\bu^\e,\rho^\e)(t)\|_{H^2(\Omega)}^2+\|(\rho^\e,(\rho^\e)^{-1})(t)\|_{L^\infty(\Omega)}^2\\
	&+\|\e(q^\e_t,\bu^\e_t)(t)\|_{H^1(\Omega)}^2+\|\e^2(q^\e_{tt},\bu^\e_{tt})(t)\|_{L^2(\Omega)}^2\\
	&+\|\e(\n^2(\mD_b\bu^\e),\n\crl\crl \bu^\e)(t)\|_{L^2(\Omega)}^2+\|\e^2\n^3 q^\e(t)\|_{L^2(\Omega)}^2,
\end{align*}
    and
\begin{align*}
	D^\e(t):=&\|(\e q^\e,\bu^\e)(t)\|_{H^3(\Omega)}^2+\|\e(q^\e_t,\bu^\e_t)(t)\|_{H^2(\Omega)}^2+\|\e^2\bu^\e_{tt}(t)\|_{H^1(\Omega)}^2\\
	&+\|\e^2(\n^3(\mD_b\bu^\e),\n^2\crl\crl \bu^\e)(t)\|_{L^2(\Omega)}^2.
\end{align*}
    The initial weighted energy is defined by
\begin{align*}
	M_0^\e:=M^\e(t)|_{t=0}.
\end{align*}
    \noindent\textbf{Remark 1.2.} The initial data are ``ill-prepared" in the sense that
\begin{align*}
	\|\dv\bu_0^\e\|_{H^1(\Omega)}+\|(\e q_t^\e(0),\e \bu_t^\e(0))\|_{H^1(\Omega)}\leq C,
\end{align*}
    under the assumption $M_0^\e\leq C$. This leads to rapidly oscillating acoustic waves in the singular limit.

    The local existence and uniqueness of the strong solution to \eqref{ts_1}-\eqref{TS_bc} for a fixed $\e$ follow from standard fixed-point arguments (cf. \cite{Hof12, MaN83, OuY22}), and the proof is omitted here.
\begin{theorem}[Local existence]\label{local}
	Let $\Omega\subset \mathbb{R}^3$ be a simply connected and bounded domain with smooth boundary $\p\Omega$. For any fixed $\e\in(0,1]$, assume that the initial datum $(q^\e_0,\bu^\e_0)$ satisfies $\rho_0^\e=b+\e \mA(b)q_0^\e>0$,
\begin{align*}
	(q^\e_0,\bu^\e_0)\in H^3(\Omega),\quad (q^\e_t(0),\bu^\e_t(0))\in H^1(\Omega),\quad (q^\e_{tt}(0),\bu^\e_{tt}(0))\in L^2(\Omega),
\end{align*}
	and that the following compatibility conditions hold:
\begin{align*}
	\bu^\e_0\cdot \bn=\bu^\e_t(0)\cdot \bn=0,\quad \crl \bu^\e_0\times \bn=\crl \bu^\e_t(0)\times \bn=0\enspace {\rm on} \enspace\p\Omega.
\end{align*}
	Then there exists a positive constant $T^*=T^*(q^\e_0,\bu^\e_0,\e)$ such that the initial boundary value problem \eqref{ts_1}-\eqref{TS_bc} admits a unique strong solution $(q^\e,\bu^\e)$ on $\Omega\times(0,T^*]$, satisfying
\begin{gather*}
	q^\e \in C([0,T^*],H^2(\Omega))\cap L^\infty(0,T^*;H^3(\Omega)),\quad \bu^\e \in C([0,T^*],H^3(\Omega))\cap L^2(0,T^*;H^4(\Omega)),\\
	q^\e_t\in C([0,T^*],H^1(\Omega))\cap L^2(0,T^*;H^2(\Omega)),\quad \bu^\e_t\in C([0,T^*],H^1(\Omega))\cap L^2(0,T^*;H^2(\Omega)),\\
	q^\e_{tt}\in L^\infty(0,T^*;L^2(\Omega)),\quad \bu^\e_{tt}\in L^\infty(0,T^*;L^2(\Omega))\cap L^2(0,T^*;H^1(\Omega)).
\end{gather*}
\end{theorem}
    \noindent\textbf{Remark 1.3.} The initial temporal derivative $q^\e_t(0)$ in Theorem \ref{local} is obtained by evaluating equation \eqref{ts_1} at $t=0$, which means that
\begin{align*}
	q^\e_t(0)=-\mA(b)^{-1}\dv(\mA(b)q_0^\e\bu_0^\e)-\frac{1}{\e}\mA(b)^{-1}\dv(b\bu_0^\e).
\end{align*}
    Moreover, $\bu^\e_t(0)$ is defined analogously by
\begin{align*}
	\bu^\e_t(0)=(\rho_0^\e)^{-1}\left[\mu\Delta\bu_0^\e+(\mu+\lambda)\n\dv\bu_0^\e-\frac{1}{\e}b\n q_0^\e-\n\pi_0^\e\right]-\bu_0^\e\cdot\n\bu_0^\e,
\end{align*}
    where
\begin{align*}
	\pi_0^\e=\frac{(\rho_0^\e)^\gamma-b^\gamma-\e bq_0^\e}{\e^2}.
\end{align*}
    The higher-order temporal derivatives $q^\e_{tt}(0)$ and $\bu^\e_{tt}(0)$ are then defined recursively by differentiating \eqref{ts_1} and \eqref{ts_2} with respect to time variable and subsequently setting $t=0$.
\begin{prop}\label{bound}
	Let $\e\in (0,1]$, and let $T^\e>0$ denote the maximal time of existence for the smooth solution $(q^\e,\bu^\e)$ to the initial-boundary value problem \eqref{ts_1}-\eqref{TS_bc}, as given in Theorem \ref{local}. Then for any $\e\in (0,1]$ and $t\in [0,T^\e)$, the following estimate holds:
\begin{align*}
	M^\e(t)\leq C_0(M_0^\e){\rm exp}\{(t^\frac{1}{4}+\e^2)C_1(M^\e(t))\}.
\end{align*}
	Here $C_0$ and $C_1$ are positive, nondecreasing, continuous functions independent of $\e$.
\end{prop}
\begin{theorem}[Uniform estimates for ill-prepared initial data]\label{estimate}
	Let $\Omega\subset\mathbb{R}^3$ be a simply connected and bounded domain with smooth boundary $\p\Omega$ and $\e\in (0,1]$ be an arbitrary constant. Let $(q^\e,\bu^\e)$ be the strong solution obtained in Theorem \ref{local} for each $\e\in(0,1]$. Suppose that $b$ satisfies \eqref{brho} and $M^\e_0\leq D_0$ for some constant $D_0>0$ independent of $\e$. Then there exist positive constants $\e_0$, $T$, and $C$, independent of $\e$, such that the following uniform estimate holds:
\begin{align}\label{unes}
	M^\e(t)\leq C,\quad \e\in(0,\e_0],\enspace t\in [0,T].
\end{align}
\end{theorem}
    Motivated by \cite{FNS08}, we introduce the weighted Helmholtz decomposition associated with the hydrostatic density $b$. For any vector $v\in H^1(\Omega)$, we write
\begin{align}\label{defH}
	v=H_b[v]+H_b^\perp[v], \quad H_{b}^\perp[v]=b\n\psi,
\end{align}
    where $\psi$ is the unique zero-mean solution to the weighted Neumann problem
\begin{align}\label{neu}
	\dv(b\n\psi)=\dv v\enspace {\rm in}\enspace \Omega,\quad b\n\psi\cdot \bn|_{\p\Omega}=v\cdot\bn,\quad \om\psi dx=0.
\end{align}
    In particular, if $v\cdot\bn|_{\p\Omega}=0$, then $\psi$ satisfies the homogeneous Neumann condition $\n\psi\cdot\bn|_{\p\Omega}=0$. We recall some standard properties of the weighted Helmholtz decomposition. With respect to the weighted inner product
\begin{align*}
	(u,v)_{1/b}:=\om\frac{u(x)\cdot v(x)}{b}dx,
\end{align*}    
    $H_b$ is the orthogonal projection onto the space of divergence-free vector fields with zero normal trace. In particular,
\begin{align*}
	\dv H_{b}[v]=0,\quad H_{b}[v]\cdot\bn|_{\p\Omega}=0,\quad H_{b}[b\n f]=0,
\end{align*}
    for every sufficiently regular scalar function $f$. In addition,
\begin{align*}
	H_{b}[v]=v
\end{align*}
    whenever $\dv v=0$ in $\Omega$ and $v\cdot\bn|_{\p\Omega}=0$. Furthermore, for $0\leq s\leq 3$,
\begin{align}\label{Hbd}
	\|H_b[v]\|_{H^s}+\|H_b^\perp[v]\|_{H^s}
	\le C\|v\|_{H^s},
\end{align}
    where $C$ depends only on $\Omega$, $s$, and the bounds in \eqref{brho}. This boundedness is a consequence of the standard elliptic regularity for the weighted Neumann problem \eqref{neu} (see Section 6.3.3 in \cite{FN17}). The uniform bounds established in \eqref{unes} yield the following anelastic limit.
    
\begin{theorem}[Singular limit]\label{limit}
	Assume that the conditions in Theorem \ref{estimate} are satisfied, and there exists a function $\mathbf{w}_0\in H^2(\Omega)$ such that as $\e\rightarrow 0$,
\begin{align*}
	\bu_0^\e \rightharpoonup \mathbf{w}_0 \enspace {\rm in} \enspace H^2(\Omega)\quad {\rm and}\quad \bu_0^\e \rightarrow \mathbf{w}_0 \enspace {\rm in} \enspace H^1(\Omega).
\end{align*}
	Then, there exists a subsequence of $\{\bu^\e\}$, still denoted by $\{\bu^\e\}$, and a function $\bv\in L^\infty(0,T;H^2(\Omega))\cap L^2(0,T;H^3(\Omega))$, such that
\begin{gather*}
	\bu^\e\overset{\ast}{\rightharpoonup} \bv\enspace {\rm in}\enspace L^\infty(0,T;H^2(\Omega)),\quad
	\bu^\e\rightharpoonup \bv\enspace {\rm in}\enspace L^2(0,T;H^3(\Omega)),\\
	H_{b}[\rho^\e \bu^\e]\rightarrow b\bv \enspace {\rm in}\enspace C([0,T],H^{2-\delta}(\Omega))\cap L^2(0,T;H^{3-\delta}(\Omega)) ,\enspace 0<\delta<1,\\
	H_{b}^\perp[\rho^\e \bu^\e]\overset{\ast}{\rightharpoonup} 0\enspace {\rm in}\enspace L^\infty(0,T;H^2(\Omega)),\quad
	H_{b}^\perp[\rho^\e \bu^\e]\rightharpoonup 0\enspace {\rm in}\enspace L^2(0,T;H^3(\Omega)).
\end{gather*}
	Moreover, there exists a scalar pressure $\Pi$, such that $(\bv,\Pi)$ is the unique strong solution to the anelastic system
\begin{gather*}
	\dv(b\bv)=0\quad {\rm in}\enspace (0,T]\times\Omega,\\
	(b\bv)_t+\dv(b\bv\otimes\bv)+b\n\Pi=\dv \mathbb{S}(\bv)
	\quad {\rm in}\enspace (0,T]\times\Omega,
\end{gather*}
	with
\begin{align*}
	\bv\cdot\bn=0,\quad \bn\times\crl\bv=0
	\quad {\rm on}\enspace (0,T]\times\p\Omega,
	\qquad
	\bv(0,x)=b^{-1}H_{b}[b\mathbf{w}_0].
\end{align*}
\end{theorem}
\begin{theorem}[Convergence rate of the projected momentum]\label{conv}
	Under the assumptions of Theorem \ref{limit}, let $\bv$ be the limiting
	strong solution on $[0,T]$. Then there exists a constant $C$, depending only on $\Omega$, $\gamma$, $\mu$, $\lambda$, $b$ and $D_0$, but independent
	of $\e$, such that
\begin{align}\label{rate}
	\|H_b[\rho^\e\bu^\e]-b\bv\|_{L_t^\infty(L^2)}+\|H_b[\rho^\e\bu^\e]-b\bv\|_{L_t^2(H^1)}\leq C\left(\|H_b[\rho_0^\e\bu_0^\e]-b\bv_0\|_{L^2}+\e\right).
\end{align}
	In particular, if $\|H_b[\rho_0^\e\bu_0^\e]-b\bv_0\|_{L^2}\leq C\e$, the left-hand side of \eqref{rate} is $O(\e)$.
\end{theorem}
    

\subsection{Difficulties and strategies} 
    The uniform estimates must account for ill-prepared data and the nonconstant hydrostatic
    density $b(x_3)$. The singular acoustic terms form the variable-coefficient pair $(b\n q^\e,\dv(b\bu^\e))$,
    rather than the constant-coefficient pair $(\n q^\e,\dv\bu^\e)$. The estimates must preserve this
    weighted structure and control the commutators generated by spatial differentiation. Two
    issues are central.
    
    \textbf{First, the classical skew-symmetric cancellation must be replaced by a weighted one.} In the standard low Mach number limit around a constant background state, the singular terms $\frac{1}{\e}\n q^\e$ and $\frac{1}{\e}\dv\bu^\e$ form a perfect skew-symmetric pair under the usual $L^2$-inner products. In the present strongly stratified regime, this cancellation is no longer available if one works directly with $\dv\bu^\e$. Indeed, the hydrostatic density $b$ varies in space, and the natural divergence variable is not $\dv\bu^\e$, but rather
\begin{align*}
	\mD_b\bu^\e=b^{-1}\dv(b\bu^\e).
\end{align*}
    Consequently, the singular terms must be tested with the weights adapted to $b$ and $\mC(b)$, where
\begin{align}
	\mC(b)=b/\mA(b)=\gamma b^{\gamma-1}.
\end{align}
    For example, the cancellation at the level of $(\n q^\e,\mD_b\bu^\e)$ is obtained from
\begin{align*}
	\frac{1}{\e}\om b\n(\mC(b)\mD_b\bu^\e)\cdot\n q^\e dx-\frac{1}{\e}\om b\n q^\e\cdot\n(\mC(b)\mD_b\bu^\e)dx=0.
\end{align*}
    These weights restore the cancellation of the singular acoustic terms. Using constant-coefficient
    test functions instead leaves terms involving derivatives of $b$ that cannot be absorbed uniformly in $\e$.


    \textbf{Second, higher-order estimates must preserve the weighted acoustic structure.} To obtain local strong solutions uniformly in $\e$, we need weighted estimates up to third
    order. A naive differentiation of the equations, especially after dividing by variable coefficients, would produce singular commutators involving derivatives of $b$, $\mC(b)$ and $\rho^\e$. These terms are not controllable by the available energy. The correct procedure is to differentiate the equations only in forms that keep the leading singular pair
\begin{align*}
	q_t^\e+\frac{1}{\e}\mC(b)\mD_b\bu^\e,\quad \bu_t^\e+\frac{1}{\e}\n q^\e
\end{align*}
    intact, and then to choose test functions carrying the same weights. In particular, we estimate $\mD_b\bu^\e$, $\mC(b)\mD_b\bu^\e$, and their derivatives rather than treating $\dv\bu^\e$ as the principal acoustic variable. The equivalence
\begin{align*}
	\dv\bu^\e=\mD_b\bu^\e-\bu^\e\cdot\n\log b
\end{align*}
    then allows us to recover the usual divergence estimates for lower-order terms.
    
\subsection{Outline of the proof of Proposition \ref{bound}}
    We estimate the acoustic variables $(\n q^\e,\mD_b\bu^\e)$ and the rotational variable $\crl\bu^\e$ separately.
    Test functions weighted by $b$ and $\mC(b)$ preserve the cancellation of the singular acoustic terms,
    while the vorticity equation contains no singular pressure or gravitational term. Under the vorticity-slip boundary conditions, elliptic div–curl estimates recover the full Sobolev norms
    from the divergence, vorticity, and lower-order velocity terms.
    
    Guided by these ideas, we now outline the sequence of estimates leading to the bound for $M^\e(t)$ in Proposition \ref{bound}. The proof proceeds from low-order to high-order norms.

    \textbf{Step 1. Uniform estimates for $\|(q^\e,\bu^\e)\|_{L_t^\infty(H^1)}$.} The basic $L^2$-estimate is obtained by testing the density equation \eqref{ts_1} with $q^\e$ and the momentum equation \eqref{ts_2} with $\bu^\e$. To estimate $(\n q^\e,\mD_b\bu^\e)$, we apply $\n$ to \eqref{TS_1} and test it by $b\n q^\e$, and test \eqref{TS_2} by $-b\n(\mC(b)\mD_b\bu^\e)$. The estimate for $\crl\bu^\e$ is then derived from the vorticity equation \eqref{crlu}, using the identity
\begin{align*}
	\crl(\rho^{-1}\n\rho^\gamma-\n g)=0
\end{align*}
    to remove the singularity.

    \textbf{Step 2. Reduction of the $H^2$-estimates.} The second-order estimate cannot be closed directly from spatial differentiation alone. Instead, the momentum equation is used to express $\n^2 q^\e$ in terms of weighted derivatives of $\mD_b\bu^\e$, vorticity terms, and $\e\bu_t$. Similarly, the density equation gives control of $\n \mD_b\bu^\e$ through $\e q_t^\e$ and lower-order terms. Schematically,
\begin{align*}
	\|\n^2q^\e\|_{L^2}\lesssim\|(\e\n^2\mD_b\bu^\e,\e\n\crl\crl\bu^\e)\|_{L^2}+\|\e\bu_t^\e\|_{H^1}+l.o.t.,
\end{align*}
    and
\begin{align*}
	\|\n\mD_b\bu^\e\|_{L^2}\lesssim\|\e q_t^\e\|_{H^1}+l.o.t..
\end{align*}
    Thus the $H^2$-estimate for $(q^\e,\bu^\e)$ is reduced to weighted estimates for the first-order temporal derivatives and to the vorticity estimates.

    \textbf{Step 3. Weighted estimates for $\|(\e q_t^\e,\e\bu_t^\e)\|_{L_t^\infty(H^1)}$.}  
    Differentiating \eqref{TS_1}-\eqref{TS_2} once in time gives the same singular acoustic pair for $(q_t^\e,\bu_t^\e)$. The $L^2$-estimate follows by the standard energy method. The $H^1$-estimate is more delicate. We test $\n\eqref{TSt}_1$ by $\e^2b\n q_t^\e$ and $\eqref{TSt}_2$ by $-\e^2b\n(\mC(b)\mD_b\bu^\e_t)$. This choice avoids dividing by $\mC(b)$, keeps the coefficient of $q_{tt}^\e$ constant, and preserves the weighted skew-symmetric cancellation. The rotational part $\e\crl\bu_t^\e$ is controlled by applying $\p_t$ to the vorticity equation \eqref{crlu}.

    \textbf{Step 4. Weighted estimates for $\|(\e^2 q^\e_{tt},\e^2\bu^\e_{tt})\|_{L_t^\infty(L^2)}$ and $\|\e^2\bu^\e_{tt}\|_{L_t^2(H^1)}$.} The derivation of these estimates is similar to that of the $L^2$-estimate.

    \textbf{Step 5. Weighted estimates for the highest-order derivatives.} The highest-order part is divided into the rotational estimate and the acoustic-divergence estimate. For the rotational component, it is more efficient to use the elliptic form of the vorticity equation rather than differentiating the vorticity equation directly. This yields
\begin{align*}
	\|\e\crl\crl\bu^\e\|_{H^1}^2+\int_0^t\|\e^2\crl\crl\bu^\e\|_{H^2}^2ds\leq C.
\end{align*}
    For the acoustic-divergence part, one first takes the inner product of $\p_i\n\eqref{TS_2}$ with $-\e^4b\p_i\n^2(\mC(b)\mD_b\bu^\e)$, and $\p_i\n^2\eqref{TS_1}$ with $\e^4b\p_i\n^2q^\e$ in $L^2$, respectively, and sums up the results. Then, the leading singular terms are canceled out, which gives the bounds of $(\e^2\n^3q^\e,\e^2\n^3\mD_b\bu^\e)$. One then tests $\p_i\n\eqref{TS_1}$ with $\e^2b\p_i\n q_t^\e$ and $\p_t\p_i\eqref{TS_2}$ with $-\e^2b\p_i\n(\mC(b)\mD_b\bu^\e)$, which yields the remaining estimates for $(\e\n^2\mD_b\bu^\e,\e\n^2q_t^\e)$ by eliminating the singular terms. Finally, using
\begin{align*}
	\dv\bu^\e=\mD_b\bu^\e-\bu^\e\cdot\n\log b
\end{align*}
    and the standard div-curl elliptic estimate under the vorticity-slip boundary condition, we recover the uniform bound of $\bu^\e$ in $L^2(0,T;H^3(\Omega))$. This closes the uniform estimate for $M^\e(t)$.

\vskip2mm
\noindent
\begin{minipage}{\linewidth}
	\centering
\begin{tikzpicture}[
	node distance=1cm and 0.8cm,
	>=Latex,
	every node/.style={align=center}
	]
	\node (n1-4) {$\|\e\n^2 q^\e_t\|_{L_t^2(L^2)}$};
	\node (n2-4) [right=of n1-4] {$\|(\e^2q^\e_{tt},\e^2\bu^\e_{tt})\|_{L_t^\infty(L^2)}$};
	\node (n3-4) [right=of n2-4] {$\|\e^2\bu^\e_{tt}\|_{L_t^2(H^1)}$};
	\node (n4-4) [right=of n3-4] {$\|\e \bu^\e_t\|_{L_t^2(H^2)}$};
	\node (n5-4) [right=of n4-4] {$\|\e \crl\crl\bu^\e_t\|_{L_t^2(L^2)}$};
	\node (n1-3) [above=of n1-4] {$\|\e\n^2\mD_b\bu^\e\|_{L_t^\infty(L^2)}$};
	\node (n2-3) [above=of n2-4] {$\|\e\n\crl\crl\bu^\e\|_{L_t^\infty(L^2)}$};
	\node (n3-3) [above=of n3-4] {$\boldsymbol{\|\e\bu^\e_t\|_{L_t^\infty(H^1)}}$};
	\node (n4-3) [above=of n4-4] {$\boldsymbol{\|\e q^\e_t\|_{L_t^\infty(H^1)}}$};
	\node (n1-2) [above=of n1-3] {$\|(q^\e,\bu^\e)\|_{L_t^\infty(H^1)}$};
	\node (n2-2) [above=of n2-3] {$\boldsymbol{\|\n^2q^\e\|_{L_t^\infty(L^2)}}$};
	\node (n3-2) [above=of n3-3] {$\|\n\dv\bu^\e\|_{L_t^\infty(L^2)}$}; 
	\node (n4-2) [above=of n4-3] {$\|\crl\crl\bu^\e\|_{L_t^\infty(L^2)}$};
	\node (n1-1) [above=of $(n1-2)!0.5!(n4-2)$] {$\|(q^\e,\bu^\e)\|_{L_t^\infty(H^2)}$};
	\draw[->,>=stealth] (n1-2) -- (n1-1);
	\draw[->,>=stealth] (n2-2) -- (n1-1);
	\draw[->,>=stealth] (n3-2) -- (n1-1);
	\draw[->,>=stealth] (n4-2) -- (n1-1);
	\draw[->,>=stealth] (n1-3) -- (n2-2)node[near start,above]{\eqref{TS_2}};
	\draw[->,>=stealth] (n2-3) -- (n2-2);
	\draw[->,>=stealth] (n3-3) -- (n2-2)node[near start,above]{\eqref{TS_2}};
	\draw[->,>=stealth] (n4-3) -- (n3-2)node[near start,above]{\eqref{TS_1}};
	\draw[<->,>=stealth,dashed] (n1-4) -- (n1-3);
	\draw[->,>=stealth] (n2-4) -- (n3-4);
	\draw[->,>=stealth] (n3-4) -- (n3-3);
	\draw[->,>=stealth] (n4-4) -- (n4-3);
	\draw[->,>=stealth] (n5-4) -- (n4-4);
\end{tikzpicture}
	\captionof{figure}{Dependencies for  the estimates of $\|(q^\e,\bu^\e)\|_{L_t^\infty(H^2)}$}
\end{minipage}

\noindent
\begin{minipage}{\linewidth}
\centering
\begin{tikzpicture}[
	node distance=1cm and 0.8cm,
	>=Latex,
	every node/.style={align=center}
	]
	\node (n1-4) {$\|\e^2\n^3\mD_b\bu^\e\|_{L_t^2(L^2)}$};
	\node (n2-4) [right=of n1-4] {$\|\e^2\crl\crl\bu^\e\|_{L_t^2(H^2)}$};
	\node (n3-4) [right=of n2-4] {$\|\e \mD_b\bu^\e\|_{L_t^\infty(H^2)}$};
	\node (n1-3) [above=of n1-4] {$\|\e^2\n^3 q^\e\|_{L_t^\infty(L^2)}$};
	\node (n2-3) [right=of n1-3] {$\|\e\n q^\e\|_{L_t^2(H^2)}$};
	\node (n3-3) [right=of n2-3] {$\|\e q^\e_t\|_{L_t^2(H^2)}$};
	\node (n1-2) [above=of n1-3] {$\|\bu^\e\|_{L_t^2(H^2)}$};
	\node (n2-2) [above=of n2-3] {$\boldsymbol{\|\dv\bu^\e\|_{L_t^2(H^2)}}$};
	\node (n3-2) [above=of n3-3] {$\|\crl\bu^\e\|_{L_t^2(H^2)}$};
	\node (n1-1) [above=of $(n1-2)!0.5!(n3-2)$] {$\|\bu^\e\|_{L_t^2(H^3)}$};
	\draw[->,>=stealth] (n1-2) -- (n1-1);
	\draw[->,>=stealth] (n2-2) -- (n1-1);
	\draw[->,>=stealth] (n3-2) -- (n1-1);
	\draw[->,>=stealth] (n3-3) -- (n2-2)node[near start,above]{\eqref{TS_1}};
	\draw[->,>=stealth] (n2-3) -- (n2-2)node[midway,left]{\eqref{TS_1}};
	\draw[<->,>=stealth,dashed] (n3-4) -- (n3-3);
	\draw[->,>=stealth] (n2-4) -- (n2-3);
	\draw[->,>=stealth] (n1-4) -- (n2-3);
	\draw[<->,>=stealth] (n1-3) -- (n1-4);
\end{tikzpicture}
	\captionof{figure}{Dependencies for  the estimates of $\|\bu^\e\|_{L_t^2(H^3)}$}
\end{minipage}

    \noindent These figures illustrate the logical structure of the uniform energy estimates. In these diagrams, single arrows denote control relationships, where the quantity at the tail controls the quantity at the head. Bidirectional arrows signify the pairs of quantities that must be estimated simultaneously to achieve the closure of estimates. Dashed arrows represent common components shared between the two figures, while boldface terms highlight more complex and challenging components of the analysis.

    The rest of this paper is organized as follows. In Section 2, we collect several auxiliary lemmas. In Section 3, we establish uniform estimates with respect to $\e$, which completes the proof of Proposition \ref{bound}. In Section 4, we complete the proofs of Theorems \ref{estimate} and \ref{limit}.

\section{Preliminaries} 
    In this section, we collect several auxiliary lemmas and inequalities used throughout the paper.
\begin{lemma}[{\cite[Part 1, Theorem 10.1]{Fdm08}}]
	Let $\Omega\subset \mathbb{R}^N$ be a bounded domain with a $C^k$-boundary, and let $\bu\in W^{k,r}(\Omega)\cap L^q(\Omega)$ with $1\leq r,q\leq \infty$. For any integer $0\leq j<k$ and any number $j/k\leq a\leq 1$, set
\begin{align*}
	\frac{1}{p}=\frac{j}{N}+a\left(\frac{1}{r}-\frac{k}{N}\right)+(1-a)\frac{1}{q}.
\end{align*}
	If $k-j-\frac{N}{r}$ is not a nonnegative integer, then
\begin{align}\label{inter}
	\|D^j\bu\|_{L^p(\Omega)}\leq C\|\bu\|_{W^{k,r}(\Omega)}^a\|\bu\|_{L^q(\Omega)}^{1-a}.
\end{align}
	If $k-j-\frac{N}{r}$ is a nonnegative integer, then \eqref{inter} holds only for $a=j/k$. The constant $C$ depends only on $\Omega$, $r$, $q$, $k$, $j$, and $a$.
\end{lemma}
    We record the following consequences of \eqref{inter} in dimension three, which will be used repeatedly:
\begin{gather*}
	\|\bu\3\leq C\|\bu\|_{H^1}^\frac{1}{2}\|\bu\2^\frac{1}{2},\quad \|\n^2 \bu\2\leq C\|\bu\|_{H^3}^\frac{1}{2}\|\n \bu\2^\frac{1}{2},\\
	\|\n^2 \bu\3\leq C\|\bu\|_{H^3}^\frac{1}{2}\|\bu\|_{H^2}^\frac{1}{2},\quad \|\n^2 \bu\4\leq C\|\bu\|_{H^3}^\frac{3}{4}\|\bu\|_{H^2}^\frac{1}{4},\\
	\|\n \bu\0\leq C\|\bu\|_{W^{2,4}}\leq C\|\bu\|_{H^3}^\frac{3}{4}\|\bu\|_{H^2}^\frac{1}{4}.
\end{gather*}

\begin{lemma}[(Moser-type composition estimate)]\label{moser}
	Let $\Omega\subset\mathbb{R}^3$ be a bounded smooth domain, $m\geq1$ be an integer, $I\subset\mathbb{R}$ be an open interval, and $K\Subset I$ be a compact interval. Let $F=F(x,z)$ be smooth on $\bar{\Omega}\times I$. If $u\in H^m(\Omega)\cap L^\infty(\Omega)$ and $u(x)\in K$ for a.e. $x\in\Omega$, then
\begin{align*}
	\|F(\cdot,u)\|_{H^m(\Omega)}
	\leq C_1\left(\|u\|_{L^\infty(\Omega)}\right)
	\left(1+\|u\|_{H^m(\Omega)}\right),
\end{align*}
	where $C_1$ is a positive nondecreasing function. More precisely, for any fixed compact interval $K'\Subset I$ with $K\Subset K'$, $C_1$ depends only on $\Omega$, $m$, $K'$, and the suprema of $\partial_x^\beta\partial_z^\ell F$ on $\bar{\Omega}\times K'$ with $|\beta|+\ell\leq m$. The notation $C_1(\|u\|_{L^\infty(\Omega)})$ records the dependence on the range of $u$.
\end{lemma}
    \noindent\textit{Proof.} We reduce the estimate to the Moser composition estimate in the whole space. Let $E$ be a standard extension operator such that, for $\tilde u=Eu$,
\begin{align*}
	\|\tilde u\|_{H^m(\mathbb{R}^3)}
	+\|\tilde u\|_{L^\infty(\mathbb{R}^3)}
	\leq C\left(\|u\|_{H^m(\Omega)}
	+\|u\|_{L^\infty(\Omega)}\right).
\end{align*}
    Setting $M=\|u\|_{L^\infty(\Omega)}$, then $\|\tilde u\|_{L^\infty(\mathbb{R}^3)}\leq CM$. Choose a compact interval $K'\Subset I$ such that $K\Subset K'$. By the smooth extension theorem in the $x$-variable and a cut-off in the $z$-variable supported in $K'$, there exists a function $\tilde F\in C^m(\mathbb{R}^3\times\mathbb{R})$ such that $\tilde F=F$ on $\bar{\Omega}\times K$ and
\begin{align*}
	\|\tilde F\|_{C^m(\mathbb{R}^3\times\mathbb{R})}
	\leq C_1(M),
\end{align*}
    where $C_1$ has the dependence stated in the lemma. Since $\tilde u=u$ on $\Omega$ and $u(x)\in K$ a.e. in $\Omega$, we have $\tilde F(x,\tilde u(x))=F(x,u(x))$ for a.e. $x\in\Omega$. Applying the Moser composition estimate \cite[Section 3.1, (14)]{BdLV11} to $\tilde F$ and $\tilde u$, we obtain
\begin{align*}
	\|F(\cdot,u)\|_{H^m(\Omega)}
	&\leq \|\tilde F(\cdot,\tilde u)\|_{H^m(\mathbb{R}^3)}\\
	&\leq C\|\tilde F\|_{C^m(\mathbb{R}^3\times\mathbb{R})}
	\left(1+\|\tilde u\|_{H^m(\mathbb{R}^3)}\right)\\
	&\leq C_1(M)\left(1+\|u\|_{H^m(\Omega)}\right).
\end{align*}
    This proves the lemma. \qed

\begin{lemma}[{\cite[Lemma 5]{Bre74}}]\label{dv-cl}
	Let $\Omega\subset \mathbb{R}^3$ be a bounded domain with smooth boundary $\p\Omega$ and its unit outward normal $\bn$. Then there exists a constant $C>0$, depending only on $\Omega$, such that for any $\bu\in H^s(\Omega)^3$ with $s\geq 1$,
\begin{align*}
	\|\bu\|_{H^s(\Omega)}\leq C\left(\|\dv \bu\|_{H^{s-1}(\Omega)}+\|\crl \bu\|_{H^{s-1}(\Omega)}+\|\bu\cdot \bn\|_{H^{s-\frac{1}{2}}(\p\Omega)}+\|\bu\|_{H^{s-1}(\Omega)}\right).
\end{align*}
\end{lemma}

\begin{lemma}[{\cite[Theorem 2.4]{Xyl07}}]\label{dv-cl'}
	Let $\Omega\subset \mathbb{R}^3$ be a bounded domain with smooth boundary $\p\Omega$ and its unit outward normal $\bn$. Then there exists a constant $C>0$, depending only on $\Omega$, such that for any $\bu\in H^s(\Omega)^3$ with $s\geq 1$,
\begin{align*}
	\|\bu\|_{H^s(\Omega)}\leq C\left(\|\dv \bu\|_{H^{s-1}(\Omega)}+\|\crl \bu\|_{H^{s-1}(\Omega)}+\|\bu\times \bn\|_{H^{s-\frac{1}{2}}(\p\Omega)}+\|\bu\|_{H^{s-1}(\Omega)}\right).
\end{align*}
\end{lemma}


\begin{lemma}[{\cite[Lemma 4]{Sim90}}]\label{cpt}
	Let $X_0\subset X\subset Y$ be three Banach spaces. Assume that $X_0$ is compactly embedded in $X$ and that $X$ is continuously embedded in $Y$. Let
\begin{align*}
	W=\{\bu\in L^p(0,T;X_0):\bu_t\in L^q(0,T;Y)\},\quad 1\leq p,q\leq +\infty.
\end{align*}
	{\rm (i)} If $p<+\infty$, then the embedding of $W$ into $L^p(0,T;X)$ is compact.\\
	{\rm (ii)} If $p=+\infty$ and $q>1$, then the embedding of $W$ into $C([0,T],X)$ is compact.
\end{lemma}

    Throughout the rest of the paper, $C$ denotes a generic positive constant independent of $\e$, while $C_0$ and $C_1$ denote positive, nondecreasing, continuous functions of their arguments. We repeatedly
    use the algebra property of $H^2(\Omega)$ in three dimensions, the embeddings $H^2(\Omega)\hookrightarrow L^\infty(\Omega)$ and $H^1(\Omega)\hookrightarrow L^6(\Omega)$, and the bounds in \eqref{brho}. 

\section{Uniform energy estimates}
    Throughout this section, we drop the superscript $\e$ on the solutions for notational brevity. We will establish that for any $\e\in(0,1]$ and $t\in[0,T]$ with $T={\rm min}\{T^\e,1\}$, the estimate
\begin{align*}
	M(t)\leq C_0(M_0){\rm exp}\{(t^\frac{1}{4}+\e^2)C_1(M(t))\}
\end{align*}
    holds. For the acoustic estimates, we write the perturbation system \eqref{ts_1}-\eqref{ts_2} as
\begin{align}\label{TS_1}
	&q_t+\frac{1}{\e}\mC(b)\mD_{b}\bu=F_q,\\ \label{TS_2}
	&\bu_t+\frac{1}{\e}\n q-\mathcal{V}=\mathcal{N},
\end{align}
    where
\begin{align}\label{F}
	F_q:&=-\mA(b)^{-1}\dv(\mA(b)q\bu),\\
	\label{V}
	\mathcal{V}:&=\rho^{-1}[\mu\Delta\bu+(\mu+\lambda)\n\dv\bu],\\ \label{N}
	\mathcal{N}:&=-\bu\cdot\n\bu-\rho^{-1}\n\pi+\mA(b)\rho^{-1}q\n q.
\end{align}

\subsection{Estimates for $\rho$ and $\pi$}
\begin{lemma}\label{RHO}
	Assume that the assumptions in Theorem \ref{local} are satisfied and that $(q,\bu)$ is sufficiently smooth with $\rho=(b^\gamma+\e b q+\e^2\pi)^\frac{1}{\gamma}$. Then for $\e\in(0,1]$ and $t\in[0,T]$, we have
\begin{gather}\label{rho}
	\|(\rho,\rho^{-1})(t)\|_{L^\infty}^2+\|\rho(t)\|_{H^2}^2\leq M_0{\rm exp}\{Ct^\frac{1}{2}M(t)\}.
\end{gather}
\end{lemma}
    \noindent\textit{Proof.} The proof is standard and we reproduce it for the convenience of readers (see \cite{Val83}). Multiplying equation \eqref{OS_1} by $\rho^{-k}$ for each integer $k\geq 2$ yields
\begin{gather*}
	\frac{1}{k-1}\om(\p_t\rho^{1-k}+\bu\cdot\n\rho^{1-k})dx-\om\rho^{1-k}\dv\bu dx=0,
\end{gather*}
    which in turn implies
\begin{gather*}
 	\frac{d}{dt}\|\rho^{-1}\|_{L^{k-1}}\lesssim\|\dv\bu\0\|\rho^{-1}\|_{L^{k-1}}.
\end{gather*}
    Applying Gr\"onwall's inequality and letting $k\rightarrow +\infty$, we obtain
\begin{align}\label{3.2}
	\|\rho^{-1}(t)\0\leq&\|\rho_0^{-1}\0{\rm exp}\left(C\int_0^t\|\bu\|_{H^3}ds\right) \nonumber\\
	\leq&M_0{\rm exp}\{Ct^\frac{1}{2}M(t)\},\quad 0\leq t\leq T.
\end{align} 
    Similarly, we can get
\begin{gather}\label{3.3}
	\|\rho(t)\0\leq\|\rho_0\0{\rm exp}\left(C\int_0^t\|\bu\|_{H^3}ds\right),\quad 0\leq t\leq T.
\end{gather}
    The $H^2$-bound of $\rho$ follows from a standard energy estimate. For any multi-index $\alpha$ with $0\leq|\alpha|\leq2$, applying $\p^\alpha$ to \eqref{OS_1}, then multiplying the result by $\p^\alpha\rho$, and integrating over $\Omega$ yield
\begin{gather*}
	\frac{d}{dt}\|\rho(t)\|_{H^2}^2\leq C\|\bu\|_{H^3}\|\rho\|_{H^2}^2,
\end{gather*}
    which implies that, for any $0<\e\leq 1$ and $0<t\leq T$,
\begin{gather}\label{3.5}
	\|\rho(t)\|_{H^2}\leq\|\rho_0\|_{H^2}{\rm exp}\left(C\int_0^t\|\bu\|_{H^3}ds\right)\leq M_0{\rm exp}\{Ct^\frac{1}{2}M(t)\}.
\end{gather}
    Collecting the estimates \eqref{3.2}-\eqref{3.5} directly yields the desired bound \eqref{rho}.\qed
    
\begin{lemma}\label{PI}
	Assume that the assumptions in Lemma \ref{RHO} are satisfied. Set
\begin{align}\label{defG}
	G=\frac{1}{\e}\left(\pi-\frac{\gamma(\gamma-1)}{2}b^{\gamma-2}\phi^2\right),
\end{align}
	where $\phi$ is defined in $\eqref{rhop}_1$. Then, for $\e\in(0,1]$ and $t\in[0,T]$, we have
\begin{gather}\label{pi1}
	\|\pi(t)\|_{H^2}+\e\|\pi_t(t)\|_{H^1}+\e^2\|\pi_{tt}\2\leq C_1(E(t)),\\\label{pi2}
	\e\|\pi(t)\|_{H^3}+\e\|\pi_t(t)\|_{H^2}\leq C_1(E(t))D^\frac{1}{2}(t),\\\label{G}
	\|G\|_{H^1}\leq C_1(E(t)).
\end{gather}
\end{lemma}
    \noindent\textit{Proof.} By Taylor's formula with integral remainder, we obtain
\begin{align}\label{pi}
	\pi=&\frac{(b+\e\phi)^\gamma-b^\gamma-\gamma\e b^{\gamma-1}\phi}{\e^2} \nonumber\\
	=&\gamma(\gamma-1)\phi^2\int_0^1(1-s)(b+s\e\phi)^{\gamma-2}ds.
\end{align}
    Define
\begin{align}\label{Bxz}
	B(x,z)=\gamma(\gamma-1)b^{\gamma-2}(x)\int_0^1(1-s)\left(1+\frac{sz}{b(x)}\right)^{\gamma-2}ds,
\end{align}
    then
\begin{align}\label{p1}
	\pi=B(x,\e\phi)\phi^2.
\end{align}
    For every $s\in[0,1]$, we have
\begin{align}\label{b}
	b+s\e\phi=(1-s)b+s\rho\geq c'>0,
\end{align}
    where $c'>0$ is determined by $c_*$ and the bound for $\rho^{-1}$ in \eqref{rho}. In particular, the graph
    $\{(x,\e\phi(x)):x\in\bar{\Omega}\}$ is separated from the singular set
    $\{(x,z): b(x)+sz=0\ {\rm for\ some}\ s\in[0,1]\}$.
    Hence, after multiplying by a smooth cut-off function in the $(x,z)$-variables and extending outside a small neighborhood
    of this graph, we may regard $B$ as a smooth function on $\bar{\Omega}\times I$, where $I$ is an open
    interval containing the range of $\e\phi$. The extended function agrees with the original one at
    $z=\e\phi(x)$, and thus Lemma \ref{moser} gives
\begin{align}\label{p2}
	\|B(x,\e\phi)\|_{H^2}\leq
	C_1(\|\e\phi\|_{L^\infty})
	(1+\|\e\phi\|_{H^2}).
\end{align}
    By H\"older's inequality, we then have
\begin{align*}
	\|\phi^2\|_{H^2}\leq C\|\phi\|_{H^2}^2.
\end{align*}
    Therefore, from \eqref{p1} we get
\begin{align}\label{pi_2}
	\|\pi\|_{H^2}=&\|B(x,\e\phi)\phi^2\|_{H^2}\nonumber\\ \leq&\|\phi\|_{H^2}^2C_1(\|\e\phi\|_{L^\infty})(1+\|\e\phi\|_{H^2}) \nonumber\\
	\leq&\|q\|_{H^2}^2C_1(\|\e q\|_{H^2}) \nonumber\\
	\leq&C_1(E(t)).
\end{align}
    Moreover, the same argument for higher-order derivatives gives
\begin{align}\label{p3}
	\|B(x,\e\phi)\|_{H^3}
	\leq C_1(E(t))(1+\e\|\phi\|_{H^3}).
\end{align}
    Due to \eqref{defq}, we have
\begin{align}\label{phiq}
	\e\|\phi\|_{H^3}\leq C\e\|q\|_{H^3}\leq CD^\frac{1}{2}(t).
\end{align}
    By using the product estimate
\begin{align*}
	\|fg\|_{H^3}\leq C(\|f\|_{H^3}\|g\|_{H^2}+\|f\|_{H^2}\|g\|_{H^3}),
\end{align*}
    \eqref{p2}, and \eqref{p3}-\eqref{phiq}, we obtain
\begin{align}
	\e\|\pi\|_{H^3}\leq&
	C\e\|B(x,\e\phi)\|_{H^3}\|\phi^2\|_{H^2}
	+C\e\|B(x,\e\phi)\|_{H^2}\|\phi^2\|_{H^3}\nonumber\\
	\leq&
	C_1(E(t))\left(\e\|\phi\|_{H^2}^2(1+\e\|\phi\|_{H^3})
	+\e\|\phi\|_{H^2}\|\phi\|_{H^3}\right)\nonumber\\
	\leq&C_1(E(t))D^\frac{1}{2}(t).
\end{align}
    Differentiating \eqref{pi} with respect to $t$ gives
\begin{align*}
	\pi_t=A_1(x,\e q)qq_t,\quad \pi_{tt}=A_1(x,\e q)qq_{tt}+A_2(x,\e q)q_t^2,
\end{align*}
    where $A_1$ and $A_2$ depend smoothly on $x$ and $\e q$, and the Sobolev norms used below are bounded by $C_1(E(t))$. Thus, we have
\begin{align*}
	&\e\|\pi_t\|_{H^1}\leq C_1(E(t))\|q\|_{H^2}\|\e q_t\|_{H^1}\leq C_1(E(t)),\\
	&\e\|\pi_t\|_{H^2}\leq C_1(E(t))\|q\|_{H^2}\|\e q_t\|_{H^2}\leq C_1(E(t))D^\frac{1}{2}(t),\\
	&\e^2\|\pi_{tt}\2\leq C_1(E(t))(\|q\|_{H^2}\|\e^2q_{tt}\2+\|\e q_t\|_{H^1}^2)\leq C_1(E(t)).
\end{align*}

    It remains to prove the estimate for $G$. From \eqref{Bxz}, we have
\begin{align*}
	B(x,0)=\gamma(\gamma-1)b^{\gamma-2}\int_0^1(1-s)\,ds
	=\frac{\gamma(\gamma-1)}{2}b^{\gamma-2}.
\end{align*}
    Therefore, we deduce from \eqref{defG} and \eqref{p1} that 
\begin{align*}
	G=\frac{B(x,\e\phi)-B(x,0)}{\e}\phi^2.
\end{align*}
    By the mean value formula with respect to the second variable of $B$, we find
\begin{align*}
	\frac{B(x,\e\phi)-B(x,0)}{\e}=\phi\int_0^1\p_zB(x,s\e\phi)ds,
\end{align*}
    and then
\begin{align}\label{Ges}
	G=C^\e(x)\phi^3,\quad C^\e(x):=\int_0^1\partial_zB(x,s\e\phi)ds.
\end{align}
    The coefficient $C^\e$ depends smoothly on $b$ and $\e\phi$. By \eqref{b} and Lemma \ref{moser},
\begin{align*}
	\|C^\e\|_{L^\infty}+\|\n C^\e\|_{L^2}\leq C_1(E(t)).
\end{align*}
    The embedding $H^2(\Omega)\hookrightarrow L^\infty(\Omega)$ and \eqref{defq} also give
\begin{align*}
	\|\phi\|_{L^\infty}+\|\n\phi\|_{L^2}\lesssim\|q\|_{H^2}\le C_1(E(t)).
\end{align*}
    Consequently, we derive from \eqref{Ges} that
\begin{gather*}
	\|G\|_{L^2}\lesssim\|C^\e\|_{L^\infty}\|\phi\|_{L^\infty}^2
	\|\phi\|_{L^2}\leq C_1(E(t)),\\
	\|\n G\|_{L^2}\lesssim\|\n C^\e\|_{L^2}\|\phi\|_{L^\infty}^3+\|C^\e\|_{L^\infty}\|\phi\|_{L^\infty}^2\|\n\phi\|_{L^2}\leq C_1(E(t)),
\end{gather*}
    which proves \eqref{G}.
    \qed
    
\subsection{$H^2$-estimates for $(q,\bu)$}
\begin{lemma}\label{QU}
	Assume that the assumptions in Lemma \ref{RHO} are satisfied. Then for $\e\in(0,1]$ and $t\in[0,T]$, we have
\begin{gather}\label{qu}
	\|(q,\bu)(t)\2^2+\int_0^t(\|\dv \bu\2^2+\|\crl \bu\2^2)ds\leq C_0(M_0){\rm exp}\{tC_1(M(t))\}.
\end{gather}
\end{lemma}
    \noindent\textit{Proof.} We take the inner product of \eqref{ts_1} with $q$ in $L^2(\Omega)$ and integrate by parts to get
\begin{align}\label{qu_1}
	\frac{1}{2}\frac{d}{dt}\om\mA(b)q^2dx+\frac{1}{\e}\om q\dv(b\bu)dx=-\frac{1}{2}\om\dv(\mA(b)\bu)q^2dx.
\end{align}
    Recalling the identity $\Delta=\n\dv-\crl\crl$ and the boundary condition \eqref{TS_bc}, we take the $L^2$-inner product of \eqref{ts_2} with $\bu$ and integrate the results by parts to obtain
\begin{align}\label{qu_2}
	\frac{1}{2}\frac{d}{dt}\om\rho|\bu|^2dx+\frac{1}{\e}\om b\n q\cdot\bu dx+\mu\om|\crl\bu|^2dx+\kappa\om|\dv\bu|^2dx=\om\pi\dv\bu dx.
\end{align}
    Consequently, combining \eqref{qu_1} and \eqref{qu_2} and using \eqref{pi1} provide that
\begin{align*}
	&\frac{d}{dt}(\|q\2^2+\|\bu\2^2)+\|\dv\bu\2^2+\|\crl\bu\2^2\\
	\leq&\|q\|_{H^2}^2\|\bu\|_{H^2}+\|\bu\|_{H^2}\|\pi\2,
\end{align*}
    which yields \eqref{qu} after integrating in $t$.\qed
    
\begin{lemma}\label{NQU}
	Assume that the assumptions in Lemma \ref{RHO} are satisfied. Then for $\e\in(0,1]$ and $t\in[0,T]$,
\begin{align}\label{nqu}
	\|\n q(t)\2^2+\|(\dv\bu,\crl\bu)(t)\2^2+\int_0^t(\|\n\dv\bu\2^2+\|\crl\crl \bu\2^2)ds\leq C_0(M_0){\rm exp}\{tC_1(M(t))\}.
\end{align}
\end{lemma}
    \noindent\textit{Proof. Step 1. Estimates for the fast variables $(\n q,\dv\bu)$.} We take the inner product of $\n\eqref{TS_1}$ with $b\n q$ in $L^2(\Omega)$ to get
\begin{align}\label{nq}
	\frac{1}{2}\frac{d}{dt}\om b|\n q|^2dx+\frac{1}{\e}\om b\n(\mC(b)\mD_b\bu)\cdot\n qdx=\om b\n F_q\cdot\n qdx. 
\end{align}
    Taking the inner product of $\eqref{TS_2}$ and $-b\n(\mC(b)\mD_b\bu)$ in $L^2(\Omega)$ and noting that
\begin{align*}
	-\om b\bu_t\cdot\n(\mC(b)\mD_b\bu)dx=&\om\dv(b\bu_t)\mC(b)\mD_b\bu dx\\
	=&\om\dv(b\bu_t)\cdot\frac{1}{\mA(b)}\dv(b\bu)dx\\
	=&\frac{1}{2}\frac{d}{dt}\om\frac{1}{\mA(b)}(\dv(b\bu))^2dx,
\end{align*}
    we have
\begin{align}\label{nu}
	\frac{1}{2}\frac{d}{dt}\om b\mC(b)\mD_b\bu^2dx-\frac{1}{\e}\om b\n q\cdot\n(\mC(b)\mD_b\bu)dx+I_V=I_N,
\end{align}
    where
\begin{align*}
	I_V:=\om b\mathcal{V}\cdot\n(\mC(b)\mD_b\bu)dx,\quad I_N:=-\om b\mathcal{N}\cdot\n(\mC(b)\mD_b\bu)dx.
\end{align*}
    Note that
\begin{align}\label{divu}
	\dv\bu=\mD_b\bu-\bu\cdot\n\log b,
\end{align}
    which gives
\begin{align}
	\rho\mV=\mu\Delta\bu+(\mu+\lambda)\n\dv\bu=\kappa\n \mD_b\bu-\mu\crl\crl\bu-\kappa\n(\bu\cdot\n\log b).
\end{align}
    Thus we have
\begin{align}\label{nu_1}
	&I_V-\kappa\om b\rho^{-1}\mC(b)|\n\mD_b\bu|^2dx \nonumber\\
	=&\kappa\om b\rho^{-1}\mD_b\bu\n\mD_b\bu\cdot\n\mC(b)dx-\mu\om b\rho^{-1}\crl\crl\bu\cdot\n(\mC(b)\mD_b\bu)dx \nonumber\\
	&-\kappa\om b\rho^{-1}\n(\bu\cdot\n\log b)\cdot\n(\mC(b)\mD_b\bu)dx \nonumber\\
	\lesssim&\|\rho^{-1}\0(\|\n\mD_b\bu\2\|\mD_b\bu\2+\|\bu\|_{H^2}\|\n\mD_b\bu\2+\|\bu\|_{H^1}\|\mD_b\bu\|_{H^1}) \nonumber\\
	\lesssim&\|\rho^{-1}\0\|\bu\|_{H^2}^2,
\end{align}
    and
\begin{align}\label{nu_2}
	I_N\lesssim&\|\bu\|_{H^2}^2\|\n\mD_b\bu\2+\|\pi\|_{H^1}\|\n\mD_b\bu\2+\|q\|_{H^2}^2\|\n\mD_b\bu\2 \nonumber\\
	\lesssim&(\|\bu\|_{H^2}^2+\|\pi\|_{H^1}+\|q\|_{H^2}^2)\|\bu\|_{H^2}.
\end{align}
    The definition of $F_q$ in \eqref{F} gives
\begin{align}\label{nq_1}
	\left|\om b\n F_q\cdot\n qdx\right|\lesssim\|q\|_{H^2}^2\|\bu\|_{H^2}.
\end{align}
    Adding \eqref{nq} and \eqref{nu}, and using \eqref{pi1} and \eqref{nu_1}-\eqref{nq_1}, we obtain
\begin{align}\label{nf'}
	&\frac{d}{dt}\om(b|\n q|^2+b\mC(b)\mD_b\bu^2)dx+\kappa\om b\rho^{-1}\mC(b)|\n \mD_b\bu|^2dx \nonumber\\
	\leq&\|\rho^{-1}\0\|\bu\|_{H^2}^2+(\|\bu\|_{H^2}^2+\|\pi\|_{H^1}+\|q\|_{H^2}^2)\|\bu\|_{H^2}.
\end{align}
    Due to \eqref{rho} and \eqref{divu}, we can infer from \eqref{nf'} that
\begin{align}\label{nf}
	\|\n q(t)\2^2+\|\dv\bu(t)\2^2+\int_0^t\|\n\dv\bu\2^2ds\leq C_0(M_0){\rm exp}\{tC_1(M(t))\}.
\end{align}

    \textit{Step 2. Estimates for the slow variable $\crl\bu$.} To derive the estimates for the vorticity $\crl\bu$, it is convenient to first rewrite the momentum equation \eqref{OS_2} in the following form
\begin{gather}\label{TS_2'}
	\bu_t+\bu\cdot\n\bu+\frac{1}{\e^2}(\rho^{-1}\n\rho^\gamma-\n g)=\rho^{-1}(-\mu\crl\crl\bu+\kappa\n\dv\bu).
\end{gather}
    The pressure and gravitational terms are curl-free:
\begin{align}\label{crl0}
	\frac{1}{\e^2}\crl(\rho^{-1}\n\rho^\gamma-\n g)=0.
\end{align}
    Using $\crl\n=0$ and $\Delta\crl=-\crl\crl\crl$, we apply ``$\crl$" to \eqref{TS_2'} and set $w=\crl\bu$ to get
\begin{align}\label{crlu}\left\{
	\begin{aligned}
	&w_t+\bu\cdot\n w-\mu\rho^{-1}\Delta w=F_w &\quad {\rm in}\enspace (0,T]\times\Omega,\\
	&\bn\times w=0 &\quad {\rm on}\enspace (0,T]\times\p\Omega,\\
	&w(0,x)=w_0:=\crl\bu_0, &\quad {\rm in}\enspace\Omega,
    \end{aligned}\right.
\end{align}
    where the boundary condition $\eqref{crlu}_2$ comes from \eqref{TS_bc}, and
\begin{align*}
	F_w:=&w\cdot\n\bu-w\dv\bu-\mu\n(\rho^{-1})\times\crl w+\kappa\n(\rho^{-1})\times\n\dv\bu.
\end{align*}
    Noting that
\begin{align*}
	\mu\om(\n\rho^{-1}\times\crl w)\cdot w dx+\mu\om(\n\rho^{-1}\times w)\cdot\crl w dx=0,
\end{align*}
    and taking the inner product of $\eqref{crlu}_1$ and $w$ in $L^2(\Omega)$, we deduce from the boundary condition $\bn\times w|_{\p\Omega}=0$ that
\begin{align}\label{ns'}
	&\frac{1}{2}\frac{d}{dt}\om|w|^2dx+\mu\om\rho^{-1}|\crl w|^2dx \nonumber\\
	=&\om w\cdot\n\bu\cdot wdx-\frac{1}{2}\om\dv\bu|w|^2dx+\kappa\om(\n(\rho^{-1})\times\n\dv\bu)\cdot wdx \nonumber\\
	\lesssim&\|\bu\|_{H^2}^3+\|\rho^{-1}\0^2\|\rho\|_{H^2}\|\bu\|_{H^2}^2.
\end{align}
    Integrating \eqref{ns'} in $[0,t]$ yields that
\begin{align}\label{ns}
	\|\crl\bu(t)\2^2+\int_0^t\|\crl\crl\bu\2^2ds\leq C_0(M_0){\rm exp}\{tC_1(M(t))\}.
\end{align}
    Combining \eqref{nf} and \eqref{ns} immediately gives \eqref{nqu}. \qed
\begin{lemma}\label{NNQU}
	Assume that the assumptions in Lemma \ref{RHO} are satisfied. Then for $\e\in(0,1]$ and $t\in[0,T]$,
\begin{gather}\label{nnq}
	\|\n^2 q(t)\2^2\leq C\|(\e\n^2\mD_b\bu,\e\n\crl\crl \bu)(t)\2^2+C_0(M_0)\|\e\bu_t(t)\|_{H^1}^2+C_0(M_0){\rm exp}\{tC_1(M(t))\},\\\label{ndu}
	\|\n\dv\bu(t)\2^2\leq C_0(M_0)\|\e q_t(t)\|_{H^1}^2+C_0(M_0){\rm exp}\{tC_1(M(t))\},\\\label{ccu}
	\|\crl\crl\bu(t)\2^2+\int_0^t\|\Delta\crl\bu\2^2ds\leq C_0(M_0){\rm exp}\{tC_1(M(t))\}.
\end{gather}
\end{lemma}
    \noindent\textit{Proof.} Applying ``$\e\n$'' to $\eqref{TS_2}$ and using \eqref{rho}, we obtain the following $L^2(\Omega)$ estimate:
\begin{align*}
	\|\n^2 q\2^2\lesssim&\|\rho^{-1}\0\|(\e\n^2\mD_b\bu,\e\n\crl\crl \bu)\2^2+\e^2\|\rho^{-1}\0^2\|\rho\|_{H^2}\|\bu\|_{H^2}\\
	&+\|\e\bu_t\|_{H^1}^2+\e^2\|\bu\|_{H^2}^4+\e^2\|\n\pi\|_{H^1}^2+\e^2\|\rho^{-1}\0\|q\|_{H^2}^2\\
	\lesssim&\|\rho^{-1}\0\|(\e\n^2\mD_b\bu,\e\n\crl\crl\bu)\2^2+C_1(E(t)),
\end{align*}
    where we have used $\Delta\bu=\n\dv\bu-\crl\crl\bu$ and Lemma \ref{PI} for $\pi$.
    Similarly, \eqref{TS_1} and Lemma \ref{inter} imply
\begin{align*}
	\|\n(\mC(b)\mD_b\bu)\2^2\lesssim&\|\e q_t\|_{H^1}^2+\e^2\|q\|_{H^2}^2\|\bu\|_{H^2}^2\\
	\lesssim&\|\e q_t\|_{H^1}^2+\e^2C_1(E(t)),
\end{align*}
    and hence \eqref{ndu} follows. Now, recalling that $\dv\crl=0$ and $\Delta w=-\crl\crl w$, we take the inner product of $\eqref{crlu}_1$ with $-\Delta w$ in $L^2(\Omega)$ to get
\begin{align*}
	&\frac{1}{2}\frac{d}{dt}\om|\crl w|^2dx+\mu\om\rho^{-1}|\Delta w|^2dx \nonumber\\
	=&\om[\crl(\bu\cdot\n\bu)+\mu\n(\rho^{-1})\times\crl w-\kappa\n(\rho^{-1})\times\n\dv\bu]\cdot\Delta wdx \nonumber\\
	\lesssim&\om\left(|\n\bu||\n\bu|+|\rho|^{-2}|\n\rho||\n^2\bu|\right)|\Delta w|dx \nonumber\\
	\lesssim&(\|\bu\|_{H^2}^2+\|\rho^{-1}\0^2\|\rho\|_{H^2}\|\bu\|_{H^2})\|\Delta w\2 \nonumber\\
	\leq&\delta\|\Delta w\2^2+C_\delta(\|\bu\|_{H^2}^2+\|\rho^{-1}\0^2\|\rho\|_{H^2}\|\bu\|_{H^2})^2,
\end{align*}
    which yields \eqref{ccu} after choosing $\delta$ sufficiently small and integrating in $t$. \qed
    
\begin{corol}\label{COR1}
	Assume that the assumptions in Lemma \ref{RHO} are satisfied. Then for $\e\in(0,1]$ and $t\in[0,T]$,
\begin{align}\label{cor1}
	\|(q,\bu)(t)\|_{H^2}^2+\int_0^t\|\crl\bu\|_{H^2}^2ds\leq C_0(M_0)(&\|(\e q_t,\e\bu_t)(t)\|_{H^1}^2+{\rm exp}\{t^\frac{1}{4}C_1(M(t))\} \nonumber\\
	&+\|(\e\n^2\mD_b\bu,\e\n\crl\crl \bu)(t)\2^2).
\end{align}
\end{corol}
    \noindent\textit{Proof.} From \eqref{TS_1} we obtain
\begin{align}
	\|\dv\bu\|_{H^1}^2\leq&C_0(M_0)(\|\e q_t\|_{H^1}^2+\e^2\|\bu\|_{H^2}^2\|q\|_{H^2}^2) \\\label{dvu}
	\leq&C_0(M_0)\|\e q_t\|_{H^1}^2+C_0(M_0){\rm exp}\{\e^2C_1(M(t))\}.
\end{align}
    Combining Lemmas \ref{dv-cl'} and \ref{RHO}, \eqref{ns}, \eqref{ccu}, and the boundary condition $\bn\times\crl\bu|_{\p\Omega}=0$ gives
\begin{align}\label{curlu}
	\|\crl\bu\|_{H^1}^2\lesssim \|\crl\crl\bu\|_{L^2}^2+\|\crl \bu\|_{L^2}^2\leq C_0(M_0){\rm exp}\{t^\frac{1}{4}C_1(M(t))\}.
\end{align}
    By using \eqref{qu}, \eqref{nqu}, \eqref{ndu}, \eqref{ccu}, \eqref{dvu}, Lemma \ref{dv-cl} for $s=1, 2$, and Lemma \ref{dv-cl'} for $s=1$, we find
\begin{align}\label{nnu}
	\|\bu\|_{H^2}^2\lesssim&(\|\dv\bu\|_{H^1}^2+\|\crl \bu\|_{H^1}^2+\|\bu\|_{H^1}^2) \nonumber\\
	\leq&C_0(M_0)(\|\e q_t\|_{H^1}^2+{\rm exp}\{(t^\frac{1}{4}+\e^2)C_1(M(t))\}).
\end{align}
    A similar argument to \eqref{curlu} yields
\begin{align}\label{crl'}
	\int_0^t\|\crl\bu\|_{H^1}^2ds\lesssim&\int_0^t\|\crl\crl \bu\|_{L^2}^2ds+\int_0^t\|\crl\bu\|_{L^2}^2ds \nonumber\\
	\leq&C_0(M_0){\rm exp}\{tC_1(M(t))\}.
\end{align}
    From \eqref{ccu}, \eqref{crl'}, and Lemmas \ref{dv-cl}-\ref{dv-cl'}, we infer that
\begin{align}\label{nncu}
	\int_0^t\|\crl\bu\|_{H^2}^2ds\lesssim&\int_0^t(\|\crl\crl \bu\|_{H^1}^2+\|\crl\bu\|_{H^1}^2)ds \nonumber\\
	\lesssim&\int_0^t(\|\crl\crl\crl\bu\2^2+\|\crl\crl\bu\cdot \bn\|_{H^\frac{1}{2}(\p\Omega)}^2+\|\crl\bu\|_{H^1}^2)ds \nonumber\\
	\lesssim&\int_0^t(\|\Delta\crl\bu\2^2+\|\bu\|_{H^2}^2)ds \nonumber\\
	\leq&C_0(M_0){\rm exp}\{t^\frac{1}{4}C_1(M(t))\}.
\end{align}
    Finally, \eqref{qu}, \eqref{nqu}, \eqref{nnq}, \eqref{nnu} and \eqref{nncu} together establish \eqref{cor1}. \qed
    
\subsection{Weighted estimates $\mathrm{I}$}
    In this subsection, we estimate $\e(q_t,\bu_t)$ in $L^2$-norm and $H^1$-norm, together with the corresponding dissipation terms. 
\begin{lemma}\label{QUt}
	Assume that the assumptions in Lemma \ref{RHO} are satisfied. Then for $\e\in(0,1]$ and $t\in[0,T]$,
\begin{gather}\label{fst}
	\|(\e q_t,\e\bu_t)(t)\2^2+\int_0^t\|(\e\dv \bu_t,\e\crl\bu_t)\2^2ds\leq C_0(M_0){\rm exp}\{tC_1(M(t))\}.
\end{gather}
\end{lemma}
    \noindent\textit{Proof.} Put
\begin{align}\label{dus}
	\mD_{b}\bu_t=b^{-1}\dv(b\bu_t).
\end{align}
    First, we differentiate \eqref{TS_1}-\eqref{TS_2} with respect to $t$ and use \eqref{TS_ic}-\eqref{TS_bc} to obtain
\begin{align}\label{TSt}\left\{
	\begin{aligned}
	&q_{tt}+\frac{1}{\e}\mC(b)\mD_{b}\bu_t=(F_q)_t &\quad {\rm in}\enspace (0,T]\times\Omega,\\
	&\bu_{tt}+\frac{1}{\e}\n q_t-(\mathcal{V})_t=(\mathcal{N})_t&\quad {\rm in}\enspace (0,T]\times\Omega,\\
	&\bu_t\cdot\bn=0,\quad \bn\times\crl\bu_t=0 &\quad {\rm on}\enspace (0,T]\times\p\Omega,\\
	&(q_t,\bu_t)(0,x)=(q_t(0),\bu_t(0)) &\quad {\rm in}\enspace \Omega,
    \end{aligned}\right.
\end{align}
    where
\begin{align*}
	q_t(0)=-\mA(b)^{-1}\dv(\mA(b)q_0\bu_0)-\frac{1}{\e}\mA(b)^{-1}\dv(b\bu_0),
\end{align*}
    and
\begin{align*}
	\bu_t(0)=(\rho_0)^{-1}\left[\mu\Delta\bu_0+(\mu+\lambda)\n\dv\bu_0-\frac{1}{\e}b\n q_0-\n\pi_0\right]-\bu_0\cdot\n\bu_0.
\end{align*}
    From \eqref{F} and \eqref{N}, we obtain
\begin{align}\label{Ft}
	&(F_q)_t=-\bu_t\cdot\n q-\bu\cdot\n q_t-q_t\dv\bu-q\dv\bu_t-q_t\bu\cdot\n\log\mA(b)-q\bu_t\cdot\n\log \mA(b),\\\label{Nt}
	&(\mathcal{N})_t=-\bu_t\cdot\n\bu-\bu\cdot\n\bu_t-(\rho^{-1})_t\n\pi-\rho^{-1}\n\pi_t+S_t\n q+S\n q_t,
\end{align}
    where
\begin{align*}
	S=\frac{\mA(b)q}{\rho},\quad(\rho^{-1})_t=\frac{-\e\mA(b)q_t}{\rho^2},\quad S_t=\frac{\mA(b)q_t}{\rho}-\frac{\e\mA(b)^2qq_t}{\rho^2}.
\end{align*}
    Multiplying $\eqref{TSt}_1$ by $\e^2\mA(b)q_t$ and $\eqref{TSt}_2$ by $\e^2b\bu_t$ in $L^2(\Omega)$, and using $\mA(b)\mC(b)=b$ yield that
\begin{align}\label{ft_1}
	\frac{1}{2}\frac{d}{dt}&\left(\e^2\om\mA(b)q_t^2dx+\e^2\om b|\bu_t|^2dx\right)-\e^2\om b(\mathcal{V})_t\cdot\bu_tdx \nonumber\\
	&=\e^2\om\mA(b)(F_q)_tq_tdx+\e^2\om b(\mathcal{N})_t\cdot\bu_tdx.
\end{align}
    Here $\bu_t\cdot\bn|_{\p\Omega}=0$ allows integration by parts, and the singular terms cancel because
\begin{align*}
	\e\om\mA(b)\mC(b)\mD_{b}\bu_tq_tdx+\e\om b\n q_t\cdot\bu_tdx=\e\om b\mD_{b}\bu_tq_tdx-\e\om q_t\dv(b\bu_t)dx=0.
\end{align*}
The right-hand side of \eqref{ft_1} can be estimated from \eqref{Ft} and \eqref{Nt} as follows:
\begin{align}\label{ft_2}
	\e^2\left|\om \mA(b)(F_q)_tq_tdx\right|\lesssim \|\e\bu_t\|_{H^1}\|\e q_t\|_{H^1}\|q\|_{H^2}+\|\e q_t\|_{H^1}^2\|\bu\|_{H^2},
\end{align}
    and
\begin{align}\label{ft_3}
	\e^2\left|\om b(\mathcal{N})_t\cdot\bu_tdx\right|\lesssim&\|\e\bu_t\|_{H^1}^2\|\bu\|_{H^2}+\|\e\pi_t\|_{H^1}\|\e\bu_t\|_{H^1} \nonumber\\
	&+\|\e q_t\|_{H^1}\|\e\bu_t\|_{H^1}(\|\pi\|_{H^1}+\|q\|_{H^2}+\|q\|_{H^2}^2).
\end{align}
    Due to \eqref{V}, we have
\begin{align}\label{Vt}
	(\mathcal{V})_t=\rho^{-1}[\mu\Delta\bu_t+(\mu+\lambda)\n\dv\bu_t]+(\rho^{-1})_t[\mu\Delta\bu+(\mu+\lambda)\n\dv\bu],
\end{align} 
    thus
\begin{align}\label{ft_4}
	&-\e^2\om b(\mathcal{V})_t\cdot\bu_tdx-\e^2\om b\rho^{-1}[\kappa|\dv\bu_t|^2+\mu|\crl\bu_t|^2]dx \nonumber\\
	=&\e^2\om[\kappa\n(b\rho^{-1})\dv\bu_t\cdot\bu_t-\mu\n(b\rho^{-1})\times\bu_t\cdot\crl\bu_t]dx \nonumber\\
	&+\e^2\om b(\rho^{-1})_t[\mu\Delta\bu+(\mu+\lambda)\n\dv\bu]\bu_tdx \nonumber\\
	\lesssim&\|\rho^{-1}\0^2\|\rho\|_{H^2}\|\e q\|_{H^2}\|\e\bu_t\|_{H^1}^2+\e\|\rho^{-1}\0^2\|\e q_t\|_{H^1}\|\bu\|_{H^2}\|\e\bu_t\|_{H^1}.
\end{align}
    Substituting \eqref{ft_2}-\eqref{ft_4} into \eqref{ft_1}, using \eqref{pi1},  and integrating the resulting inequality over $[0,t]$, we obtain \eqref{fst}. \qed
    
\begin{lemma}\label{NQUt}
	Assume that the assumptions in Lemma \ref{RHO} are satisfied. Then for $\e\in(0,1]$ and $t\in[0,T]$,
\begin{align}\label{fsnt}
	\|\e\n q_t(t)\2^2+&\|(\e\dv\bu_t,\e\crl \bu_t)(t)\2^2+\int_0^t(\|\e\n\dv \bu_t\2^2+\|\e\crl\crl\bu_t\2^2)ds \nonumber\\
	&\leq C_0(M_0){\rm exp}\{(t^\frac{1}{4}+\e^2)C_1(M(t))\}.
\end{align}
\end{lemma}
    \noindent\textit{Proof.} 
    Multiplying $\n\eqref{TSt}_1$ by $\e^2b\n q_t$ and $\eqref{TSt}_2$ by $-\e^2b\n(\mC(b)\mD_{b}\bu_t)$ yields
\begin{align}\label{ntf_1}
	&\frac{1}{2}\frac{d}{dt}\left(\e^2\om b|\n q_t|^2dx+\e^2\om b\mC(b)\mD_{b}\bu_t^2dx\right)+\e^2\om b(\mathcal{V})_t\cdot\n(\mC(b)\mD_{b}\bu)_tdx \nonumber\\
	=&\e^2\om b\n(F_q)_t\cdot\n q_tdx-\e^2\om b(\mathcal{N})_t\cdot\n(\mC(b)\mD_{b}\bu)_tdx,
\end{align}
    where we have canceled the singular terms by observing that
\begin{align*}
	\e\int_\Omega b\n(\mC(b)\mD_b\bu_t)\cdot\n q_tdx-\e\int_\Omega b\n q_t\cdot\n(\mC(b)\mD_{b}\bu)_tdx=0,
\end{align*}
    and have used that
\begin{align*}
	-\e^2\om b\bu_{tt}\cdot\n(\mC(b)\mD_{b}\bu)_tdx=\e^2\om\dv(b\bu_{tt})\mC(b)\mD_{b}\bu_tdx=\frac{1}{2}\frac{d}{dt}\om\e^2b\mC(b)\mD_{b}\bu_t^2dx.
\end{align*}
    Using \eqref{Ft} and \eqref{Nt}, we have
\begin{align}\label{ntf_2}
	\e^2\left|\om b\n(F_q)_t\cdot\n q_tdx\right|\lesssim\|\bu\|_{H^3}\|\e q_t\|_{H^1}^2+\|q\|_{H^2}\|\e\bu_t\|_{H^1}\|\e q_t\|_{H^1},
\end{align}
    and
\begin{align}\label{ntf_3}
	\e^2\left|\om b(\mathcal{N})_t\cdot\n(\mC(b)\mD_{b}\bu)_tdx\right|\lesssim&\|\e\bu_t\|_{H^1}\|\e\bu_t\|_{H^2}\|\bu\|_{H^2}+\|\e\pi_t\|_{H^1}\|\e\bu_t\|_{H^2} \nonumber\\
	&+\|\e q_t\|_{H^1}\|\e\bu_t\|_{H^2}(\|\pi\|_{H^1}+\|q\|_{H^2}+\|q\|_{H^2}^2).
\end{align}
    Rewrite \eqref{Vt} as
\begin{align*}
	(\mathcal{V})_t=&\rho^{-1}[\kappa\n \mD_{b}\bu_t-\mu\crl\crl\bu_t-\kappa\n(\bu_t\cdot\n\log b)]\\
	&+(\rho^{-1})_t[\kappa\n \mD_{b}\bu-\mu\crl\crl\bu-\kappa\n(\bu\cdot\n\log b)].
\end{align*}
    Then we have
\begin{align}\label{ntf_4}
	&\e^2\om b(\mathcal{V})_t\cdot\n(\mC(b)\mD_{b}\bu)_tdx-\e^2\om\kappa b\mC(b)\rho^{-1}|\n\mD_{b}\bu_t|^2dx \nonumber\\
	\lesssim&\|\e\bu_t\|_{H^1}\|\e\bu_t\|_{H^2}+\e\|\e q_t\|_{H^1}\|\e\bu_t\|_{H^2}\|\bu\|_{H^3}.
\end{align}
    Substituting \eqref{ntf_2}-\eqref{ntf_4} into \eqref{ntf_1} and integrating the resulting inequality over $[0,t]$, we obtain
\begin{align}\label{fnt}
	\|\e\n q_t(t)\2^2+\|\e\mD_{b}\bu_t(t)\2^2+\int_0^t\|\e\n \mD_{b}\bu_t\2^2ds\leq C_0(M_0){\rm exp}\{(t^\frac{1}{2}+\e^2)C_1(M(t))\}.
\end{align}
    
    We next estimate $\e\crl\bu_t$ in $L^2$-norm. By \eqref{crl0},
\begin{align*}
	\crl\p_t\left(\frac{1}{\e}b\rho^{-1}\n q+\rho^{-1}\n\pi\right)=0.
\end{align*}
    Applying the operator $\crl\p_t$ to \eqref{TS_2'} yields
\begin{align}\label{st}\left\{
    \begin{aligned}
 	&w_{tt}+\mu\crl(\rho^{-1}\crl w_t)=(F_w)_t &\quad {\rm in}\enspace (0,T]\times\Omega,\\
	&\bn\times w_t=0 &\quad {\rm on}\enspace (0,T]\times\p\Omega,\\
	&w_t(0,x):=\crl\bu_t(0) &\quad {\rm in}\enspace\Omega,
	\end{aligned}\right.
\end{align}
    where
\begin{align*}
	w_t(0):=&\crl[\rho_0^{-1}(\mu\Delta\bu_0+(\mu+\lambda)\n\dv\bu_0)]-\crl(\bu_0\cdot\n\bu_0)\\
	&-\frac{1}{\e}\n\left(\frac{b}{\rho_0}\right)\times\n q_0-\n(\rho_0^{-1})\times\n\pi_0,
\end{align*}
    and
\begin{align*}
	(F_w)_t:=&-\crl(\bu_t\cdot\n\bu+\bu\cdot\n\bu_t)+\kappa\crl(\rho^{-1}\n\dv\bu_t)+\crl((\rho^{-1})_t[\mu\Delta\bu+(\mu+\lambda)\n\dv\bu]).
\end{align*}
    By virtue of \eqref{inter} and \eqref{rho}, we estimate
\begin{align}\label{st_2}
	\left|\e^2\om\crl(\mathcal{N})_t\cdot w_tdx\right|\lesssim&\|\e \bu_t\|_{H^1}^\frac{1}{2}\|\e\bu_t\|_{H^2}^\frac{1}{2}(\|\e\bu_t\|_{H^1}\|\bu\|_{H^2}+\|\e q_t\|_{H^1}\|\pi\|_{H^2}) \nonumber\\
	&+\|\e \bu_t\|_{H^1}^\frac{1}{2}\|\e\bu_t\|_{H^2}^\frac{1}{2}\|\e q_t\|_{H^1}(\|q\|_{H^2}+\|q\|_{H^2}^2).
\end{align}
    Taking the inner product of $\eqref{st}_1$ and $\e^2w_t$ in $L^2(\Omega)$ gives
\begin{align}\label{st_1}
	&\frac{1}{2}\frac{d}{dt}\om|\e w_t|^2dx+\mu\om\rho^{-1}|\e\crl w_t|^2dx \nonumber\\
	=&\e^2\om (F_w)_t\cdot w_tdx \nonumber\\
	\lesssim&\left|\e^2\om\crl(\mathcal{N})_t\cdot w_tdx\right|+\|\e\bu_t\|_{H^2}^\frac{3}{2}\|\e\bu_t\|_{H^1}^\frac{1}{2}+\|\e q_t\|_{H^1}\|\e\bu_t\|_{H^1}^\frac{1}{2}\|\e\bu_t\|_{H^2}^\frac{1}{2}\|\bu\|_{H^3}.
\end{align}
    Substituting \eqref{st_2} into \eqref{st_1} and integrating over $[0,t]$, we get
\begin{align}\label{snt}
	\|\e\crl \bu_t(t)\2^2+\int_0^t\|\e\crl\crl\bu_t\2^2ds\leq C_0(M_0){\rm exp}\{t^\frac{1}{4}C_1(M(t))\}.
\end{align}
    Thus, \eqref{fsnt} follows from \eqref{fnt} and \eqref{snt}. \qed
    
\begin{corol}\label{COR2}
	Assume that the assumptions in Lemma \ref{RHO} are satisfied. Then for $\e\in(0,1]$ and $t\in[0,T]$,
\begin{align}\label{cor2}
	\|(\e q_t,\e\bu_t)(t)\|_{H^1}^2+\int_0^t\|\e\bu_t\|_{H^2}^2ds\leq C_0(M_0){\rm exp}\{(t^\frac{1}{4}+\e^2)C_1(M(t))\}.
\end{align}
\end{corol}
    \noindent\textit{Proof.} Applying Lemma \ref{dv-cl} to $\e\bu_t$, we derive that
\begin{align}\label{cor2_1}
	\|\e\bu_t\|_{H^1}^2\leq C(\|\e\dv\bu_t\2^2+\|\e\crl\bu_t\2^2+\|\e\bu_t\2^2),
\end{align}
    and
\begin{align}\label{cor2_2}
	\int_0^t\|\e\bu_t\|_{H^2}^2ds\leq C\int_0^t(\|\e\n \mD_{b}\bu_t\2^2+\|\e\crl\crl \bu_t\2^2+\|\e\dv\bu_t\2^2+\|\e\crl\bu_t\2^2+\|\e\bu_t\2^2)ds.
\end{align}
    Hence, \eqref{cor2} follows from \eqref{fst}, \eqref{fsnt}, \eqref{cor2_1}, and \eqref{cor2_2}. \qed
    
\subsection{Weighted estimates $\mathrm{II}$}
    In this subsection, we estimate $\e^2(q_{tt},\bu_{tt})$ in $L^2$-norm. 
\begin{lemma}\label{Lem3.6}
	Assume that the assumptions in Lemma \ref{RHO} are satisfied. Then for $\e\in(0,1]$ and $t\in[0,T]$,
\begin{gather}\label{qutt}
	\|(\e^2 q_{tt},\e^2\bu_{tt})(t)\2^2+\int_0^t\|\e^2\bu_{tt}\|_{H^1}^2ds\leq C_0(M_0){\rm exp}\{t^\frac{1}{4}C_1(M(t))\}.
\end{gather}
\end{lemma}
    \noindent\textit{Proof.} 
First, we differentiate \eqref{TSt} with respect to $t$ to obtain
\begin{align}\label{TStt}\left\{
\begin{aligned}
	&q_{ttt}+\frac{1}{\e}\mC(b)\mD_{b}\bu_{tt}=(F_q)_{tt} &\quad {\rm in}\enspace (0,T]\times\Omega,\\
	&\bu_{ttt}+\frac{1}{\e}\n q_{tt}-(\mathcal{V})_{tt}=(\mathcal{N})_{tt} &\quad {\rm in}\enspace (0,T]\times\Omega,\\
	&\bu_{tt}\cdot\bn=0,\quad \bn\times\crl\bu_{tt}=0 &\quad {\rm on}\enspace (0,T]\times\p\Omega,\\
	&(q_{tt},\bu_{tt})(0,x)=(q_{tt}(0),\bu_{tt}(0)) &\quad {\rm in}\enspace\Omega,
\end{aligned}\right.
\end{align}
    where the initial quantities $q_{tt}(0)$ and $\bu_{tt}(0)$ are, respectively, defined by
\begin{align*}
	q_{tt}(0):=&-\mA(b)^{-1}\dv(\mA(b)q_t(0)\bu_0+\mA(b)q_0\bu_t(0))-\frac{1}{\e}\mA(b)^{-1}\dv(b\bu_t(0)),
\end{align*}
    and
\begin{align*}
	\bu_{tt}(0):=&\rho_0^{-1}\left[\mu\Delta\bu_t(0)+(\mu+\lambda)\n\dv\bu_t(0)-\frac{1}{\e}b\n q_t(0)-\n\pi_t(0)\right]\\
	&-\rho_0^{-1}\rho_t(0)\left[(\bu_t(0)+\bu_0\cdot\n\bu_0)\right]-\bu_t(0)\cdot\n \bu_0-\bu_0\cdot\n\bu_t(0).
\end{align*}
    The second-order time derivative of $F_q$, $\mV$, and $\mN$ are
\begin{align*}
	(F_q)_{tt}:=&-\bu_{tt}\cdot\n q-2\bu_t\cdot\n q_t-\bu\cdot\n q_{tt}-q_{tt}\dv\bu-2q_t\dv\bu_t-q\dv\bu_{tt}\\
	&-q_{tt}\bu\cdot\n\log\mA(b)-2q_t\bu_t\cdot\n\log\mA(b)-q\bu_{tt}\cdot\n\log\mA(b),
\end{align*}
\begin{align*}
	(\mathcal{V})_{tt}:=&\rho^{-1}[\mu\Delta\bu_{tt}+(\mu+\lambda)\n\dv\bu_{tt}]+2(\rho^{-1})_t[\mu\Delta\bu_t+(\mu+\lambda)\n\dv\bu_t]\nonumber\\
	&+(\rho^{-1})_{tt}[\mu\Delta\bu+(\mu+\lambda)\n\dv\bu],
\end{align*}
    and
\begin{align*}
	(\mathcal{N})_{tt}:=&-\bu_{tt}\cdot\n\bu-2\bu_t\cdot\n\bu_t-\bu\cdot\n\bu_{tt}-(\rho^{-1})_{tt}\n\pi-2(\rho^{-1})_t\n\pi_t-\rho^{-1}\n\pi_{tt}\\
	&+S_{tt}\n q+2S_t\n q_t+S\n q_{tt},
\end{align*}
    where
\begin{align*}
	&(\rho^{-1})_{tt}=-\e \mA(b)q_{tt}\rho^{-2}+2\e^2\mA(b)^2q_t^2\rho^{-3},\\ &S_{tt}=\frac{\mA(b)q_{tt}}{\rho}-\frac{2\e\mA(b)^2q_t^2+\e\mA(b)^2qq_{tt}}{\rho^2}+\frac{2\e^2\mA(b)^3qq_t^2}{\rho^3}.
\end{align*}
    
    Taking the $L^2$-inner product of $\eqref{TStt}_1$ with $\e^4\mA(b)q_{tt}$ and $\eqref{TStt}_2$ with $\e^4b\bu_{tt}$, we deduce that
\begin{align}\label{ftt_1}
	&\frac{1}{2}\frac{d}{dt}\left(\e^4\om\mA(b)q_{tt}^2dx+\e^4\om b|\bu_{tt}|^2dx\right)-\e^4\om b(\mathcal{V})_{tt}\cdot\bu_{tt}dx\nonumber\\
	=&\e^4\om\mA(b)(F_q)_{tt}q_{tt}dx+\e^4\om b(\mathcal{N})_{tt}\cdot\bu_{tt}dx,
\end{align}
    where we have canceled the singular terms by
\begin{align*}
	\e^3\om\mA(b)\mC(b)\mD_{b}\bu_{tt}q_{tt}dx+\e^3\om b\n q_{tt}\cdot\bu_{tt}dx=\e^3\om b\mD_{b}\bu_{tt}q_{tt}dx-\e^3\om q_{tt}\dv(b\bu_{tt})dx=0.
\end{align*}
    According to the definitions of $(F_q)_{tt}$ and $(\mathcal{N})_{tt}$, we have
\begin{align}\label{ftt_2}
	\e^4\left|\om\mA(b)(F_q)_{tt}q_{tt}dx\right|\lesssim&\e^4\om(|\bu_{tt}||\n q|+|\bu_t||\n q_t|+|q_{tt}||\dv\bu|+|q_t||\dv\bu_t|)|q_{tt}|dx \nonumber\\
	&+\e^4\om(|q_{tt}||\bu|+|q_t||\bu_t|+|q||\bu_{tt}|)|q_{tt}|dx \nonumber\\
	\lesssim&\|\e^2\bu_{tt}\|_{H^1}\|\e^2q_{tt}\2\|q\|_{H^2}+\|\e^2q_{tt}\2\|\e \bu_t\|_{H^2}\|\e q_t\|_{H^1}+\|\e^2q_{tt}\2^2\|\bu\|_{H^3},
\end{align}
    and
\begin{align}\label{ftt_3}
	&\e^4\left|\om b(\mathcal{N})_{tt}\cdot\bu_{tt}dx\right| \nonumber\\
	\lesssim&\e^4\om(|\bu_{tt}||\n\bu|+|\bu_t||\n\bu_t|+|\dv(b\bu)||\bu_{tt}|)|\bu_{tt}|dx+\e^4\om|\pi_{tt}||\dv(b\rho^{-1}\bu_{tt})|dx \nonumber\\
	&+\e^4\om(|\e q_{tt}|+|\e q_t|^2)|\n\pi||\bu_{tt}|dx+\e^4\om|\e q_t||\n\pi_t||\bu_{tt}|dx \nonumber\\
	&+\e^4\om(|q_{tt}|+\e|q_t|^2+\e|q||q_{tt}|+\e^2|q||q_t|^2)|\n q||\bu_{tt}|dx \nonumber\\
	&+\e^4\om(|q_t|+\e|q||q_t|)|\n q_t||\bu_{tt}|dx+\e^4\om|q_{tt}||\dv(bS\bu_{tt})|dx \nonumber\\
	\lesssim&\|\e^2\bu_{tt}\|_{L^2}^2\|\bu\|_{H^3}+\|\e\bu_t\|_{H^1}^2\|\e^2\bu_{tt}\|_{H^1}+\|\e^2\bu_{tt}\|_{H^1}\|\e^2\pi_{tt}\|_{L^2} \nonumber\\
	&+\|\e^2\bu_{tt}\|_{L^2}(\|\e^2q_{tt}\|_{L^2}\|\e\pi\|_{H^3}+\|\e q_t\|_{H^1}\|\pi\|_{H^2}+\|\e q_t\|_{H^1}\|\e\pi_t\|_{H^1}) \nonumber\\
	&+\|\e^2q_{tt}\|_{L^2}\|q\|_{H^2}\|\e^2\bu_{tt}\|_{L^2}^\frac{1}{2}\|\e^2\bu_{tt}\|_{H^1}^\frac{1}{2}+\|\e q_t\|_{H^1}^2\|q\|_{H^2}\|\e^2\bu_{tt}\|_{L^2}.
\end{align}
    The principal part of the viscous term gives
\begin{align}\label{ftt_4}
	&-\e^4\om b(\mathcal{V})_{tt}\cdot\bu_{tt}dx-\e^4\om b\rho^{-1}[\kappa|\dv\bu_{tt}|^2+\mu|\crl\bu_{tt}|^2]dx \nonumber\\
	=&\e^4\om[\kappa\n(b\rho^{-1})\dv\bu_{tt}\cdot\bu_{tt}-\mu\n(b\rho^{-1})\times\bu_{tt}\cdot\crl\bu_{tt}]dx \nonumber\\
	&+2\e^4\om b(\rho^{-1})_t[\mu\Delta\bu_t+(\mu+\lambda)\n\dv\bu_t]\cdot\bu_{tt}dx \nonumber\\
	&+\e^4\om b(\rho^{-1})_{tt}[\mu\Delta\bu+(\mu+\lambda)\n\dv\bu]\cdot\bu_{tt}dx \nonumber\\
	\lesssim&\|\e^2\bu_{tt}\|_{L^2}^\frac{1}{2}\|\e^2\bu_{tt}\|_{H^1}^\frac{3}{2}\|\e q\|_{H^2}+\|\e q_t\|_{H^1}\|\e\bu_t\|_{H^2}\|\e^2\bu_{tt}\|_{L^2}^\frac{1}{2}\|\e^2\bu_{tt}\|_{H^1}^\frac{1}{2} \nonumber\\
	&+\|\e^2q_{tt}\|_{L^2}\|\e^2\bu_{tt}\|_{L^2}^\frac{1}{2}\|\e^2\bu_{tt}\|_{H^1}^\frac{1}{2}\|\bu\|_{H^3}+\|\e q_t\|_{H^1}^2\|\bu\|_{H^3}\|\e^2\bu_{tt}\|_{L^2}.
\end{align}
    Therefore, integrating the sum of \eqref{ftt_1}-\eqref{ftt_4} and applying Lemma \ref{dv-cl} yields \eqref{qutt}. \qed
    
\subsection{Weighted estimates $\mathrm{III}$}
    This subsection is devoted to the derivations of the higher order weighted estimates.
\begin{lemma}\label{NCCU}
	Assume that the assumptions in Lemma \ref{RHO} are satisfied. Then for $\e\in(0,1]$ and $t\in[0,T]$,
\begin{align}\label{nccu}
	\|\e\crl\crl\bu(t)\|_{H^1}^2+\int_0^t\|\e^2\crl\crl\bu\|_{H^2}^2ds\leq C_0(M_0){\rm exp}\{(t^\frac{1}{4}+\e^2)C_1(M(t))\}.
\end{align}
\end{lemma}
    \noindent\textit{Proof.} We use the elliptic equation for the vorticity. Equation $\eqref{crlu}_1$ gives
\begin{align}
	\mu\Delta w=\rho w_t+\rho\crl(\bu\cdot\n\bu)+\mu\rho\n(\rho^{-1})\times\crl w-\kappa\rho\n(\rho^{-1})\times\n\dv\bu.
\end{align}
    By \eqref{defr}, \eqref{rho}, and the interpolation inequality,
\begin{align*}
	\|\e\n\rho^{-1}\|_{L^\infty}
	\leq C_0(M_0)(1+\e^2\|\n q\|_{L^\infty}+\e^2\|q\|_{L^\infty})
	\leq C_0(M_0)C_1(E(t)).
\end{align*}
    By virtue of \eqref{fsnt}, \eqref{ndu}-\eqref{ccu}, Lemma \ref{PI}, and the identity $\Delta=\n\dv-\crl\crl$, we get
\begin{align*}
	\|\e\Delta\crl\bu\2^2
	\leq&C_0(M_0)C_1(E(t))[\|\e\crl\bu_t\2^2+\e^2\|\bu\|_{H^2}^4+\|\n\dv\bu\2^2+\|\crl\crl \bu\2^2\\
	&\quad+\e^2\|q\|_{H^2}^4+\|\pi\|_{H^1}^2]\\
	\leq&C_0(M_0){\rm exp}\{(t^\frac{1}{4}+\e^2)C_1(M(t))\},
\end{align*}
    which gives that
\begin{align}\label{nccu_1}
	\|\e\crl\crl\bu\|_{H^1}^2\leq&C\|\e\Delta\crl\bu\2^2+C\|\e\crl\crl\bu\cdot \bn\|_{H^\frac{1}{2}(\p\Omega)}^2+C\|\e\crl\crl\bu\2^2 \nonumber\\
	\leq&\|\e\Delta\crl\bu\2^2+\e^2\|\bu\|_{H^2}^2 \nonumber\\
	\leq&C_0(M_0){\rm exp}\{(t^\frac{1}{4}+\e^2)C_1(M(t))\}.
\end{align}
    Then we take the $L^2$-inner product of $\e^2\n\eqref{crlu}_1$ and $\e^2\n\Delta\crl\bu$ to obtain
\begin{align}\label{n3cu}
	&\mu\int_0^t\|\e^2\sqrt{\rho^{-1}}\n\Delta\crl\bu\2^2ds \nonumber\\
	=&\int_0^t\om\e^2[\n w_t+\n(\bu\cdot\n w)+\mu\rho^{-2}\n\rho\cdot\Delta w-\n F_w]\cdot\e^2\n\Delta\crl\bu dxds \nonumber\\
	\leq&C_0(M_0)\int_0^t\e\|\e^2\n\Delta\crl\bu\2[\|\e\crl \bu_t\|_{H^1}+\e\|\bu\|_{H^2}\|\bu\|_{H^3}+(\|\e^2q\|_{H^3}+\|q\|_{H^2}^2)\|\bu\|_{H^3}]ds \nonumber\\
	&+C_0(M_0)\int_0^t\e\|\e^2\n\Delta\crl\bu\2[\|q\|_{H^2}^2+\e\|\pi\|_{H^2}^2]ds \nonumber\\
	\leq&\delta\int_0^t\|\e^2\n\Delta\crl\bu\2^2ds+\e^2C_\delta C_1(E(t))D(t),
\end{align}
    which leads to
\begin{align}
	\int_0^t\|\e^2\n\Delta\crl\bu\2^2ds\leq C_0(M_0){\rm exp}\{(t^\frac{1}{4}+\e^2)C_1(M(t))\}.
\end{align}
    By \eqref{n3cu}, Lemma \ref{dv-cl}, and the identity $\dv\crl=0$,
\begin{align}\label{nccu_2}
	\int_0^t\|\e^2\crl\crl\bu\|_{H^2}^2ds\leq&C\int_0^t\left(\|\e^2\Delta\crl \bu\|_{H^1}^2+\e^2\|\crl\crl\bu\cdot\bn\|_{H^\frac{3}{2}(\p\Omega)}^2+\e^2\|\crl\crl \bu\|_{H^1}^2\right)ds \nonumber\\
	\leq&C\int_0^t\left(\|\e^2\n\Delta\crl\bu\2^2+\e^2\|\bu\|_{H^3}^2\right)ds \nonumber\\
	\leq&C_0(M_0){\rm exp}\{(t^\frac{1}{4}+\e^2)C_1(M(t))\}.
\end{align}
    Combining \eqref{nccu_1} and \eqref{nccu_2} completes the proof of \eqref{nccu}. \qed
    
\begin{lemma}\label{NNDU}
	Assume that the assumptions in Lemma \ref{RHO} are satisfied. Then for $\e\in(0,1]$ and $t\in[0,T]$,
\begin{align}\label{spt}
	\|(\e^2\n^3 q,\e\n^2\mD_{b}\bu)(t)\2^2+\int_0^t\|(\e^2\n^3\mD_{b}\bu,\e\n^2 q_t,\e\n^3 q)\2^2ds\leq C_0(M_0){\rm exp}\{(t^\frac{1}{4}+\e^2)C_1(M(t))\}.
\end{align}
\end{lemma}
    \noindent\textit{Proof.} 
    \textit{Step 1. Estimate for $\e^2\n^3q$ and $\e^2\n^3\mD_{b}\bu$}. Testing $\p_i\n\eqref{TS_2}$ with $-\e^4b\p_i\n^2(\mC(b)\mD_{b}\bu)$ and $\p_i\n^2\eqref{TS_1}$ with $\e^4b\p_i\n^2 q$, and using the cancellation
\begin{align*}
	\frac{1}{\e}\om b\p_i\n^2(\mC(b)\mD_{b}\bu)\cdot\p_i\n^2q dx-\frac{1}{\e}\om b\p_i\n^2q\cdot\p_i\n^2(\mC(b)\mD_{b}\bu)dx=0,
\end{align*}
    we obtain
\begin{align}\label{nqh_1}
	&\frac{1}{2}\frac{d}{dt}\om\e^4b|\p_i\n^2 q|^2dx+\e^4\om b\p_i\n\mathcal{V}\cdot\p_i\n^2(\mC(b)\mD_{b}\bu)dx \nonumber\\
	=&\e^4\om b\p_i\n\bu_t\cdot\p_i\n^2(\mC(b)\mD_{b}\bu)dx+\e^4\om b\p_i\n^2 F_q\cdot\p_i\n^2 qdx-\e^4\om b\p_i\n\mathcal{N}\cdot\p_i\n^2(\mC(b)\mD_{b}\bu)dx.
\end{align}
    Moreover, putting $\beta=\n\log b$, the identity
\begin{align*}
	\mu\Delta\bu+(\mu+\lambda)\n\dv\bu=\kappa\n \mD_{b}\bu-\mu\crl\crl\bu-\kappa\n(\bu\cdot\beta)
\end{align*}
    implies
\begin{align*}
	&\e^4\om b\p_i\n\mathcal{V}\cdot\p_i\n^2(\mC(b)\mD_{b}\bu)dx\\
	=&\kappa\e^4\om b\p_i\n(\rho^{-1}\n\mD_{b}\bu)\cdot\p_i\n^2(\mC(b)\mD_{b}\bu)dx-\mu\e^4\om b\p_i\n(\rho^{-1}\crl\crl\bu)\cdot\p_i\n^2(\mC(b)\mD_{b}\bu)dx\\
	&-\kappa\e^4\om b\p_i\n(\rho^{-1}\n (\bu\cdot\beta))\cdot\p_i\n^2(\mC(b)\mD_{b}\bu)dx\\
	=&:J_1+J_2+J_3,
\end{align*}
    where
\begin{align}\label{nqh_2}
	&J_1-\kappa\e^4\om b\mC(b)\rho^{-1}|\p_i\n^2\mD_{b}\bu|^2dx \nonumber\\
	\lesssim&\e^4\om b[\p_i\n,\rho^{-1}]\n\mD_{b}\bu\cdot\p_i\n^2(\mC(b)\mD_{b}\bu)dx+\e^4\om b\p_i\n(\rho^{-1}\n\mD_{b}\bu)\cdot[\p_i\n^2,\mC(b)]\mD_{b}\bu dx \nonumber\\
	\lesssim&\|\e^2\n^3\mD_{b}\bu\|_{L^2}(\e\|\e^2q\|_{H^3}\|\n\mD_{b}\bu\|_{L^2}^\frac{1}{2}\|\n\mD_{b}\bu\|_{H^1}^\frac{1}{2}+\|\rho\|_{H^2}\|\n^2\mD_{b}\bu\|_{L^2}^\frac{1}{2}\|\e^2\n^2\mD_{b}\bu\|_{H^1}^\frac{1}{2}) \nonumber\\
	&+\e^2\|\rho\|_{H^2}\|\e^2\n\mD_{b}\bu\|_{H^2}\|\mD_{b}\bu\|_{H^2} \nonumber\\
	\leq&C_0(M_0)C_1(E(t))(D^\frac{1}{2}(t)+\e^2D(t)),
\end{align}
    and
\begin{align}\label{nqh_3}
	J_2+J_3\lesssim&\|\rho\|_{H^2} \|\e^2\n^3\mD_{b}\bu\|_{L^2}\|\e^2\crl\crl\bu\|_{H^2}+\e^2\|\rho\|_{H^2}\|\mD_{b}\bu\|_{H^2}\|\e^2\crl\crl\bu\|_{H^2} \nonumber\\
	&+\e^2\|\rho\|_{H^2}\|\bu\|_{H^3}\|\e^2\mD_{b}\bu\|_{H^3} \nonumber\\
	\leq&C_0(M_0)C_1(E(t))(D^\frac{1}{2}(t)+\e^2D(t)).
\end{align}
    Using 
\begin{align*}
	F_q=-\bu\cdot\n q-q\mD_{b}\bu+(\gamma-1)q\bu\cdot\beta,
\end{align*}
    we have
\begin{align}\label{nqh_4}
	&\left|\e^4\om b\p_i\n^2 F_q\cdot\p_i\n^2 qdx\right| \nonumber\\
	\lesssim&\left|\e^4\om\dv(b\bu)|\p_i\n^2q|^2dx\right|+\e^4\|[\p_i\n^2,\bu\cdot\n]q\|_{L^2}\|\p_i\n^2q\|_{L^2} \nonumber\\
	&+\e^4\|\p_i\n^2(q\mD_{b}\bu)\|_{L^2}\|\p_i\n^2q\|_{L^2}+\e^4\|\p_i\n^2(q\bu\cdot\beta)\|_{L^2}\|\p_i\n^2q\|_{L^2} \nonumber\\
	\lesssim&\|\bu\|_{H^3}\|\e^2\n^3q\|_{L^2}^2+\|\bu\|_{H^2}^\frac{1}{2}\|\bu\|_{H^3}^\frac{1}{2}\|\e^2q\|_{H^3}\|\e^2\n^3q\|_{L^2} \nonumber\\
	&+\|\bu\|_{H^3}\|\e^2q\|_{H^3}\|\e^2\n^3q\|_{L^2}+\|q\|_{H^2}\|\e^2\n^3\mD_{b}\bu\|_{L^2}\|\e^2\n^3q\|_{L^2} \nonumber\\
	\leq&C_0(M_0)C_1(E(t))D^\frac{1}{2}(t).
\end{align}
    By virtue of the identity
\begin{align*}
	\mathcal{N}=-\bu\cdot\n\bu-\rho^{-1}\n\pi+\mA(b)\rho^{-1}q\n q,
\end{align*}
    we obtain
\begin{align}\label{nqh_5}
	&\left|\e^4\om b\p_i\n\bu_t\cdot\p_i\n^2(\mC(b)\mD_{b}\bu)dx\right|+\left|\e^4\om b\p_i\n\mathcal{N}\cdot\p_i\n^2(\mC(b)\mD_{b}\bu)dx\right| \nonumber\\
	\lesssim&\|\e^2\mD_{b}\bu\|_{H^3}(\e\|\e\bu_t\|_{H^2}+\|\bu\|_{H^2}\|\bu\|_{H^3}+\|\rho\|_{H^2}\|\e\pi\|_{H^3}+\|\rho\|_{H^2}\|q\|_{H^2}\|\e^2 q\|_{H^3}) \nonumber\\
	\leq&C_0(M_0)C_1(E(t))(D^\frac{1}{2}(t)+\e^2D(t)).
\end{align}
    Substituting \eqref{nqh_2}-\eqref{nqh_5} into \eqref{nqh_1} and integrating in $t$, we get
\begin{align}\label{3q4u}
	\|\e^2\n^3q(t)\|_{L^2}^2+\int_0^t\|\e^2\n^3\mD_{b}\bu\|_{L^2}^2ds\leq C_0(M_0){\rm exp}\{(t^\frac{1}{4}+\e^2)C_1(M(t))\}.
\end{align}
    An argument similar to that for \eqref{nnq}, based on \eqref{TS_2}, \eqref{nccu} and \eqref{3q4u}, yields
\begin{align}\label{3q}
	\int_0^t\|\e\n^3 q\2^2ds\leq&C_0(M_0)\int_0^t(\|\e^2\n\dv \bu\|_{H^2}^2+\|\e^2\crl\crl\bu\|_{H^2}^2)ds \nonumber\\
	&+C_0(M_0)\int_0^t\e^2(\|\e \bu_t\|_{H^2}^2+\|\bu\|_{H^2}^2\|\bu\|_{H^3}^2+\|\n\pi\|_{H^2}^2)ds \nonumber\\
	\leq&C_0(M_0){\rm exp}\{(t^\frac{1}{4}+\e^2)C_1(M(t))\}.
\end{align}
    \textit{Step 2. Estimate for $\e\n^2\mD_{b}\bu$ and $\e\n^2q_t$}. Taking the $L^2$-inner product of $\p_i\n\eqref{TS_1}$ with $\e^2b\p_i\n q_t$ and $\p_t\p_i\eqref{TS_2}$ with $-\e^2b\p_i\n(\mC(b)\mD_{b}\bu)$ yields
\begin{align}\label{hq}
	&\frac{1}{2}\frac{d}{dt}\om\kappa b\mC(b)\rho^{-1}|\e\p_i\n \mD_{b}\bu|^2dx+\om b|\e\p_i\n q_t|^2dx \nonumber\\
	=&\e^2\om b\p_i\bu_{tt}\cdot\p_i\n(\mC(b)\mD_{b}\bu)dx+\e^2\om b\p_i\n F_q\cdot\p_i\n q_tdx \nonumber\\
	&-\e^2\om b\p_i(\mathcal{N})_t\cdot\p_i\n(\mC(b)\mD_{b}\bu)dx-R_\beta,
\end{align}
    where
\begin{align}
	R_\beta=\e^2\om b\p_i(\mathcal{V})_t\cdot\p_i\n(\mC(b)\mD_{b}\bu)dx-\frac{1}{2}\frac{d}{dt}\om\kappa b\mC(b)\rho^{-1}|\e\p_i\n\mD_{b}\bu|^2dx.
\end{align}
    We can obtain the estimate for $\int_0^t\|\e^2\crl\crl\bu_t\|_{H^1}ds$ directly from $\eqref{st}_1$ with the help of Lemma \ref{dv-cl}, \eqref{cor1}, \eqref{fst}, \eqref{nccu}, \eqref{cor2}, \eqref{qutt} and \eqref{3q4u} as follows:
\begin{align}\label{cccu}
	\int_0^t\|\e^2\crl\crl\bu_t\|_{H^1}ds\equiv&\int_0^t\|\e^2\crl w_t\|_{H^1}ds\nonumber\\
	\lesssim&\int_0^t\|\e^2\crl\crl w_t\|_{L^2}ds+\int_0^t\|\e^2\crl w_t\cdot\bn\|_{H^\frac{1}{2}(\p\Omega)}ds+\int_0^t\|\e^2\crl w_t\2ds\nonumber\\
	\lesssim&\int_0^t\|\e^2\crl\crl w_t\2ds+\e^2\int_0^t\|\e\bu_t\|_{H^2}ds\nonumber\\
	\leq&C_0(M_0){\rm exp}\{(t^\frac{1}{4}+\e^2)C_1(M(t))\},
\end{align}
    where $w_t=\crl\bu_t$. In light of the identity
\begin{align*}
	(\mathcal{V})_t=&\rho^{-1}[\mu\Delta\bu_t+(\mu+\lambda)\n\dv\bu_t]+(\rho^{-1})_t[\mu\Delta\bu+(\mu+\lambda)\n\dv\bu]\\
	=&\rho^{-1}[\kappa\n\mD_{b}\bu_t-\mu\crl\crl\bu_t-\kappa\n(\bu_t\cdot\beta)]+(\rho^{-1})_t[\mu\Delta\bu+(\mu+\lambda)\n\dv\bu],
\end{align*}
    we have
\begin{align}\label{hq_1}
	R_\beta=&-\frac{\kappa}{2}\om\p_t(b\mC(b)\rho^{-1})|\e\p_i\n\mD_{b}\bu|^2dx-\mu\e^2\om b\rho^{-1}\p_i\crl\crl\bu_t\cdot\p_i\n(\mC(b)\mD_{b}\bu)dx \nonumber\\
	&+\e^2\om b\p_i\rho^{-1}[\kappa\n\mD_{b}\bu_t-\mu\crl\crl\bu_t-\kappa\n(\bu_t\cdot\beta)]\cdot\p_i\n(\mC(b)\mD_{b}\bu)dx \nonumber\\
	&+\e^2\om\kappa b\rho^{-1}\p_i\n\mD_{b}\bu_t\cdot(\p_i\n(\mC(b)\mD_{b}\bu)-\mC(b)\p_i\n\mD_{b}\bu)dx \nonumber\\
	&-\kappa\e^2\om b\rho^{-1}\p_i\n(\bu_t\cdot\beta)\cdot\p_i\n(\mC(b)\mD_{b}\bu)dx \nonumber\\
	&+\e^2\om b\p_i[(\rho^{-1})_t(\kappa\n\mD_{b}\bu-\mu\crl\crl\bu-\kappa\n(\bu\cdot\beta))]\cdot\p_i\n(\mC(b)\mD_{b}\bu)dx \nonumber\\
	\lesssim&\|\e q_t\|_{H^1}\|\e\mD_{b}\bu\|_{H^2}^\frac{1}{2}\|\e\mD_{b}\bu\|_{H^3}^\frac{1}{2}+\|\e^2\crl\crl\bu_t\|_{H^1}\|\mD_{b}\bu\|_{H^2}+\e\|\e q\|_{H^3}\|\e\mD_{b}\bu_t\|_{H^1}\|\mD_{b}\bu\|_{H^2} \nonumber\\
	&+\e\|\e\mD_{b}\bu_t\|_{H^1}\|\mD_{b}\bu\|_{H^2}+\e\|\e q_t\|_{H^1}\|\bu\|_{H^3}\|\mD_{b}\bu\|_{H^2}.
\end{align}
    The other terms satisfy
\begin{align}\label{hq_2}
	\left|\e^2\om b\p_i\bu_{tt}\cdot\p_i\n(\mC(b)\mD_{b}\bu)dx\right|
	\lesssim\|\e^2\bu_{tt}\|_{H^1}\|\mD_{b}\bu\|_{H^2}
	\leq\delta\|\mD_{b}\bu\|_{H^2}^2+C_\delta\|\e^2\bu_{tt}\|_{H^1}^2,
\end{align}
    and
\begin{align}\label{hq_3}
	&\left|\e^2\om b\p_i\n F_q\cdot\p_i\n q_tdx\right| \nonumber\\
	\lesssim&\|\bu\|_{H^2}\|\e\n^3q\|_{L^2}\|\e\n^2q_t\|_{L^2}+\|\e[\p_i\n,\bu\cdot\n]q\|_{L^2}\|\e\p_i\n q_t\|_{L^2} \nonumber\\
	&+\|\e\p_i\n(q\mD_{b}\bu)\|_{L^2}\|\e\p_i\n q_t\|_{L^2}+\e^2\|\p_i\n(q\bu\cdot\beta)\|_{L^2}\|\e\p_i\n q_t\|_{L^2} \nonumber\\
	\lesssim&\|\bu\|_{H^2}\|\e\n^3q\|_{L^2}\|\e\n^2q_t\|_{L^2}+\e\|q\|_{H^2}\|\bu\|_{H^2}^\frac{1}{2}\|\bu\|_{H^3}^\frac{1}{2}\|\e\n^2q_t\|_{L^2} \nonumber\\
	&+\e\|q\|_{H^2}\|\mD_{b}\bu\|_{H^2}\|\e\n^2q_t\|_{L^2}+\e\|q\|_{H^2}\|\bu\|_{H^2}\|\e\n^2q_t\|_{L^2} \nonumber\\
	\leq&\delta\|\e\n^2q_t\|_{L^2}^2+C_\delta C_0(M_0)C_1(E(t))(D^\frac{1}{2}(t)+\e^2D(t)),
\end{align}
    Due to the definition of $\mathcal{N}$, we get
\begin{align}\label{hq_4}
	&\left|\e^2\om b\p_i(\mathcal{N})_t\cdot\p_i\n(\mC(b)\mD_{b}\bu)dx\right| \nonumber\\
	\lesssim&\|\bu\|_{H^2}\|\e\bu_t\|_{H^2}\|\e\mD_{b}\bu\|_{H^2}+\|\rho\|_{H^2}\|\e^2\pi_t\|_{H^2}\|\mD_{b}\bu\|_{H^2}+\|\rho\|_{H^2}\|\e q_t\|_{H^1}\|q\|_{H^2}\|\e\mD_{b}\bu\|_{H^2} \nonumber\\
	&+\|(\rho^{-1})_t\|_{H^1}\|\e\pi\|_{H^3}\|\e\mD_{b}\bu\|_{H^2}+\|(\rho^{-1})_t\|_{H^1}\|q\|_{H^2}\|\e q\|_{H^3}\|\e\mD_{b}\bu\|_{H^2} \nonumber\\
	\leq&C_0(M_0)C_1(E(t))(D^\frac{1}{2}(t)+\e^2D(t)).
 \end{align}
    Corollaries \ref{COR1}-\ref{COR2} together with \eqref{3q} imply
\begin{align}\label{nndu}
	\int_0^t\|\dv\bu\|_{H^2}^2ds\lesssim&\int_0^t(\|\e q_t\|_{H^2}^2+\|\bu\|_{H^2}^2\|\e\n q\|_{H^2}^2+\e^2\|q\|_{H^2}^2\|\bu\|_{H^3}^2)ds \nonumber\\
	\leq&C\int_0^t\|\e\n^2q_t\2^2ds+C_0(M_0){\rm exp}\{(t^\frac{1}{4}+\e^2)C_1(M(t))\}.
\end{align}
    Since $\mD_{b}\bu=\dv\bu+\bu\cdot\n\log b$, it follows from \eqref{nndu} and the previous estimates that
\begin{align}\label{hbd}
	\int_0^t\|\mD_{b}\bu\|_{H^2}^2ds
	\leq C\int_0^t\|\e\n^2q_t\2^2ds
	+C_0(M_0){\rm exp}\{(t^\frac{1}{4}+\e^2)C_1(M(t))\}.
\end{align}
    
    Finally, substituting \eqref{hq_1}-\eqref{hq_4} into \eqref{hq}, using \eqref{qutt}, \eqref{cccu}, and \eqref{hbd}, and choosing $\delta>0$ sufficiently small to absorb the corresponding terms into the left-hand side, we deduce that
\begin{align}\label{3u4q}
	\|\e\n^2\mD_{b}\bu(t)\2^2+\int_0^t\|\e\n^2 q_t\2^2ds\leq C_0(M_0){\rm exp}\{(t^\frac{1}{4}+\e^2)C_1(M(t))\},
\end{align}
    which, in conjunction with \eqref{3q4u} and \eqref{3q}, completes the proof of \eqref{spt}. \qed
    
\begin{corol}\label{COR3}
	Assume that the assumptions in Lemma \ref{RHO} are satisfied. Then for $\e\in(0,1]$ and $t\in[0,T]$,
\begin{align}\label{cor3}
	\int_0^t\|\bu\|_{H^3}^2ds\leq C_0(M_0){\rm exp}\{(t^\frac{1}{4}+\e^2)C_1(M(t))\}.
\end{align}
\end{corol}
    \noindent\textit{Proof.} By Lemma \ref{dv-cl}, \eqref{cor1}, \eqref{nndu}-\eqref{3u4q}, and the boundary condition $\bu\cdot\bn|_{\p\Omega}=0$,
\begin{align}
	\int_0^t\|\bu\|_{H^3}^2ds\lesssim&\int_0^t\left(\|\dv \bu\|_{H^2}^2+\|\crl\bu\|_{H^2}^2+\|\bu\cdot \bn\|_{H^\frac{5}{2}(\p\Omega)}^2+\|\bu\|_{H^2}^2\right)ds \nonumber\\
	\leq&C_0(M_0){\rm exp}\{(t^\frac{1}{4}+\e^2)C_1(M(t))\}.
\end{align}
    This finishes the proof of Corollary \ref{COR3}. \qed
    
    Combining the estimates of this section proves Proposition \ref{bound}. 
    
\section{Proofs of the main theorems}
\subsection{Proof of Theorem \ref{estimate}}
    This argument follows the approaches developed by Alazard \cite{Ala06}, and M\'etivier and Schochet \cite{Met01}. Let $T^\e<+\infty$ be the maximal existence time of the solution from Theorem \ref{local}, then for all $0\leq t\leq {\rm min}\{T^\e,1\}$, the estimate
\begin{align}\label{1}
	M^\e(t)\leq C_0(M_0^\e){\rm exp}\{(t^\frac{1}{4}+\e^2)C_1(M^\e(t))\}
\end{align}
    holds, with $M_0^\e\leq D_0$ for $0<\e\leq 1$. We now fix a constant $D$ satisfying $C_0(D_0)<D<+\infty$, and choose $\e_1\in(0,1]$ and $T_1\in(0,1]$ such that
\begin{align}\label{2}
	C_0(D_0){\rm exp}\{(T_1^\frac{1}{4}+\e_1^2)C_1(D)\}<D.
\end{align}
    From inequality \eqref{1}, \eqref{2} and the initial condition $M^\e(0)=M_0^\e$, we find that $M^\e(t)\neq D$ holds for any $t<{\rm min}\{T^\e,T_1\}$. Indeed, if $M^\e(t)=D$, \eqref{1} and \eqref{2} imply
\begin{align*}
	D=M^\e(t)\leq C_0(M_0^\e){\rm exp}\{(t^\frac{1}{4}+\e^2)C_1(D)\}<D.
\end{align*}
    
    Since $M^\e(0)\leq D_0<D$, the continuity of $M^\e(t)$ gives
\begin{align}
	M^\e(t)<D \quad{\rm for}\enspace t<{\rm min}\{T^\e,T_1\},\enspace 0<\e\leq \e_1.
\end{align}
    A standard continuation argument then shows that $T^\e\geq T_1$ for any $0<\e\leq \e_1$. Hence, $M^\e(t)<D$ holds for any $t\in [0,T_1]$, with $T_1$ independent of $\e\in(0,\e_1]$.
    
\subsection{Proof of Theorem \ref{limit}}
    Let $\{(q^\e,\bu^\e)\}_{\e>0}$ be the family of local strong solutions to \eqref{TS_1}-\eqref{TS_2} with initial data satisfying the assumptions of Theorem \ref{estimate}. The uniform estimates in \eqref{unes} yield a subsequence, still denoted by $\e$, along which the following convergences holds as $\epsilon \to 0$:
\begin{gather}\label{ulmt}
	\bu^\e\rightharpoonup\bv\enspace {\rm in}\enspace L^2(0,T;H^3(\Omega)),\quad \bu^\e\overset{\ast}{\rightharpoonup}\bv\enspace{\rm in}\enspace L^\infty(0,T;H^2(\Omega)),\\
	q^\e\overset{\ast}{\rightharpoonup} q\enspace {\rm in}\enspace L^\infty(0,T;H^2(\Omega)),\\
	p^\e=b^\gamma+\e bq^\e+\e^2\pi^\e\rightarrow b^\gamma\enspace {\rm in}\enspace L^\infty(0,T;H^2(\Omega)),
\end{gather}
    and therefore
\begin{align}\label{rlmt}
	\rho^\e=(b^\gamma+\e bq^\e+\e^2\pi^\e)^\frac{1}{\gamma}\rightarrow b\enspace {\rm in}\enspace L^\infty(0,T;H^2(\Omega)).
\end{align}
    It is useful to rewrite the momentum equation as:
\begin{align}\label{eq1}
	\p_t(\rho^\e\bu^\e)+\dv(\rho^\e\bu^\e\otimes \bu^\e)+\frac{1}{\e}b\n q^\e+\n\pi^\e=\mu\Delta\bu^\e+(\mu+\lambda)\n\dv\bu^\e.
\end{align}
    
    By \eqref{unes} and \eqref{Hbd}, $H_{b}[\rho^\e\bu^\e]$ is uniformly bounded in $L^\infty(0,T;H^2(\Omega))\cap L^2(0,T;H^3(\Omega))$. Applying $H_b$ to \eqref{eq1} and following the decomposition in \cite{Mas07}, we get
\begin{align}\label{Peq}
	\p_t(H_{b}[\rho^\e\bu^\e])+&H_b\left[\dv\left(\frac{H_b[\rho^\e\bu^\e]\otimes H_b[\rho^\e\bu^\e]}{b}+\frac{H_b[\rho^\e\bu^\e]\otimes H^\perp_{b}[\rho^\e\bu^\e] }{b}+\frac{H^\perp_{b}[\rho^\e\bu^\e]\otimes H_b[\rho^\e\bu^\e]}{b}\right)\right] \nonumber\\
	+&H_b\left[\dv\left(\frac{H^\perp_{b}[\rho^\e\bu^\e]\otimes H^\perp_{b}[\rho^\e\bu^\e]}{b}\right)+\n\pi^\e\right] \nonumber\\
	=&H_{b}[\mu\Delta\bu^\e+(\mu+\lambda)\n\dv\bu^\e]-H_b\left[\dv\left(\left(\frac{1}{\rho^\e}-\frac{1}{b}\right)\rho^\e\bu^\e\otimes\rho^\e\bu^\e\right)\right].
\end{align}
    From \eqref{rhop} and \eqref{defq}, we obtain
\begin{align*}
	\left\|\frac{1}{\rho^\e}-\frac{1}{b}\right\|_{L_t^\infty(H^2)}\leq C\e,
\end{align*}
    which, combined with \eqref{Hbd}, yields
\begin{align}\label{err1}
	\left\|H_b\left[\dv\left(\left(\frac{1}{\rho^\e}-\frac{1}{b}\right)\rho^\e\bu^\e\otimes\rho^\e\bu^\e\right)\right]\right\|_{L_t^\infty(H^1)}\leq C\e.
\end{align}
    An application of Lemma \ref{cpt} then implies that, up to a subsequence, as $\e\rightarrow 0$,
\begin{align}
	H_{b}[\rho^\e\bu^\e]\rightarrow H_{b}[b \bv]\enspace {\rm in}\enspace L^2(0,T;H^{3-\delta}(\Omega))\cap C([0,T],H^{2-\delta}(\Omega)),\quad \forall \delta\in(0,1).
\end{align}
    Passing to the limit in the continuity equation by \eqref{ulmt} and \eqref{rlmt} gives $\dv(\rho^\e\bu^\e)\rightarrow \dv(b\bv)$ and $-\rho^\e_t\rightarrow 0$ in $\mD'((0,T)\times\Omega)$ as $\e\rightarrow 0$. By uniqueness of the limit, we obtain
\begin{align}
	\dv(b\bv)=0\quad {\rm a.e.}\enspace {\rm in}\enspace (0,T)\times\Omega.
\end{align}
    Together with the boundary condition $\bv\cdot\bn|_{\p\Omega}=0$, this yields $H_{b}[b\bv]=b\bv$ and $H^\perp_{b}[b\bv]=0$. Therefore, as $\e\rightarrow 0$,
\begin{align}\label{Pu}
	H_{b}[\rho^\e\bu^\e]\rightarrow b\bv\enspace{\rm in}\enspace L^2(0,T;H^{3-\delta}(\Omega))\cap C([0,T],H^{2-\delta}(\Omega)),\quad \forall \delta\in(0,1),
\end{align}
    and
\begin{align}\label{Qu}
	H^\perp_{b}[\rho^\e\bu^\e]\rightharpoonup 0\enspace {\rm in}\enspace L^2(0,T;H^{3}(\Omega)),\quad H^\perp_{b}[\rho^\e\bu^\e]\overset{\ast}{\rightharpoonup} 0\enspace {\rm in}\enspace L^\infty(0,T;H^{2}(\Omega)).
\end{align}
    
    We decompose the convection tensor as
\begin{align}\label{lmt_1}
	\frac{\rho^\e\bu^\e\otimes\rho^\e \bu^\e}{b}=&\frac{H_{b}[\rho^\e\bu^\e]\otimes H_{b}[\rho^\e\bu^\e]}{b}+\frac{H_{b}[\rho^\e\bu^\e]\otimes H^\perp_{b}[\rho^\e\bu^\e]}{b} \nonumber\\
	&+\frac{H^\perp_{b}[\rho^\e\bu^\e]\otimes H_{b}[\rho^\e\bu^\e]}{b}+\frac{H^\perp_{b}[\rho^\e\bu^\e]\otimes H^\perp_{b}[\rho^\e\bu^\e]}{b}.
\end{align}
     By \eqref{Pu}, the uniform $H^2$-bounds, and the Sobolev product estimate, we obtain
\begin{align}\label{lmt_2}
	\frac{H_{b}[\rho^\e\bu^\e]\otimes H_{b}[\rho^\e\bu^\e]}{b}\rightarrow b\bv\otimes\bv \enspace {\rm in}\enspace C([0,T],H^{2-\delta}(\Omega)).
\end{align}
    The strong convergence in \eqref{Pu} and the weak convergences in \eqref{Qu} imply that as $\e\rightarrow 0$,
\begin{align}\label{lmt_3}
	\frac{H_{b}[\rho^\e\bu^\e]\otimes H^\perp_{b}[\rho^\e\bu^\e]}{b}+\frac{H^\perp_{b}[\rho^\e\bu^\e]\otimes H_{b}[\rho^\e\bu^\e]}{b}\rightharpoonup 0&\enspace {\rm in}\enspace L^2(0,T;H^{3}(\Omega)), \nonumber\\
	\frac{H_{b}[\rho^\e\bu^\e]\otimes H^\perp_{b}[\rho^\e\bu^\e]}{b}+\frac{H^\perp_{b}[\rho^\e\bu^\e]\otimes H_{b}[\rho^\e\bu^\e]}{b} \overset{\ast}{\rightharpoonup} 0&\enspace {\rm in}\enspace L^\infty(0,T;H^{2}(\Omega)).
\end{align}
    
    For the acoustic component, write $H_b^\perp[\rho^\e\bu^\e]=b\n\psi^\e$. Since $\dv(H_b[\rho^\e\bu^\e])=0$, the continuity equation \eqref{ts_1} gives
\begin{align}\label{ac_1}
	\p_t\phi^\e+\frac{1}{\e}\dv(b\n\psi^\e)=0.
\end{align}
    Projecting the momentum equation \eqref{eq1} onto the acoustic component yields
\begin{align}\label{ac_2}
	\p_t(b\n\psi^\e)+\frac{\gamma}{\e}b\n(b^{\gamma-2}\phi^\e)=bF^\e,
\end{align}
    where 
\begin{align*}
	bF^\e:=H_b^\perp\left[\mu\Delta\bu^\e+(\mu+\lambda)\n\dv\bu^\e-\dv\left(\frac{\rho^\e\bu^\e\otimes\rho^\e\bu^\e}{\rho^\e}\right)-\n\pi^\e\right]
\end{align*}
    is uniformly bounded in $L^\infty(0,T;L^2(\Omega))\cap L^2(0,T;H^1(\Omega))$. The weak convergence of the acoustic component does not determine its quadratic self-interaction. Inspired by Section 3 in \cite{Mas07}, we combine this term with the pressure correction $\n\pi^\e$ and define
\begin{align}\label{Res}
 	\mathcal{R}^\e=&\dv(b\n\psi^\e\otimes \nabla\psi^\e) +\n\pi^\e \nonumber\\
	=&b\n\left(\frac{|\n\psi^\e|^2}{2}\right)+\dv(b\n\psi^\e)\n\psi^\e+\n\pi^\e \nonumber\\
	=&b\n\left(\frac{|\n\psi^\e|^2}{2}\right)+\dv(b\n\psi^\e)\n\psi^\e+\frac{\gamma(\gamma-1)}{2}\n(b^{\gamma-2}(\phi^\e)^2)+\e\n G^\e.
\end{align}
    Substituting \eqref{ac_1} into \eqref{Res} gives
\begin{align}\label{R_1}
	\mathcal{R}^\e=b\n\left(\frac{|\n\psi^\e|^2}{2}\right)-\e\p_t\phi^\e\n\psi^\e+\frac{\gamma(\gamma-1)}{2}\n(b^{\gamma-2}(\phi^\e)^2)+\e\n G^\e.
\end{align}
    The product rule gives
\begin{align*}
	-\e(\p_t\phi^\e)\n\psi^\e=-\e\p_t(\phi^\e\n\psi^\e)+\e\phi^\e\p_t(\n\psi^\e).
\end{align*}
    Using \eqref{ac_2}, we obtain
\begin{align}\label{R_2}
	-\e(\p_t\phi^\e)\n\psi^\e=&-\e\p_t(\phi^\e\n\psi^\e)+\e\phi^\e\left(-\frac{\gamma}{\e}\n(b^{\gamma-2}\phi^\e)+F^\e\right) \nonumber\\
	=&-\e\p_t(\phi^\e\n\psi^\e)-\gamma\phi^\e\n (b^{\gamma-2}\phi^\e)+\e \phi^\e F^\e.
\end{align}
    Notice that
\begin{align}\label{R_3}
	-\gamma\phi^\e\n(b^{\gamma-2}\phi^\e)=\frac{\gamma(\gamma-2)}{2}b\n(b^{\gamma-3}(\phi^\e)^2)-\frac{\gamma(\gamma-1)}{2}\n(b^{\gamma-2}(\phi^\e)^2).
\end{align}
    Putting \eqref{R_2} and \eqref{R_3} into \eqref{R_1}, we can rewrite the acoustic residual $\mathcal{R}^\e$ as
\begin{align}\label{lmt_4}
	\mathcal{R}^\e=&b\n\left(\frac{|\n\psi^\e|^2}{2}+\frac{\gamma(\gamma-2)}{2}b^{\gamma-3}(\phi^\e)^2\right)-\e\p_t(\phi^\e\n\psi^\e)+\e\phi^\e F^\e+\e\n G^\e.
\end{align}
    Applying $H_b$ to \eqref{lmt_4} eliminates the weighted gradient term and gives
\begin{align}\label{HbR}
	H_b[\mathcal{R}^\e]=H_b[-\e\p_t(\phi^\e\n\psi^\e)+\e\phi^\e F^\e+\e\n G^\e].
\end{align}
    
    Since $\phi^\e\n\psi^\e$ is uniformly bounded in $L^\infty(0,T;L^2(\Omega))$, integration by parts in time yields 
\begin{align*}
	H_b[\e\p_t
	(\phi^\e\n\psi^\e)]\rightarrow 0\enspace {\rm in}\enspace\mathcal D'((0,T)\times\Omega).
\end{align*}
    Moreover, the uniform bounds for $\phi^\e$ in $L^\infty(0,T;H^2(\Omega))$, $F^\e$ in $L^\infty(0,T;L^2(\Omega))$, and $G^\e$ in $L^\infty(0,T;H^1(\Omega))$, together with \eqref{Hbd}, imply
\begin{align}\label{err2}
	\|H_b\left[\e\phi^\e F^\e
	+\e\n G^\e]\right\|_{L_t^2(L^2)}\leq C\e.
\end{align}
    Consequently, \eqref{lmt_4} gives
\begin{align}\label{Hr}
	H_b[\mathcal R^\e]\rightarrow 0\enspace{\rm in}\enspace\mathcal D'((0,T)\times\Omega).
\end{align}
    Moreover, we have
\begin{align*}
	H_b[\n\pi^\e]=H_b\left[b\n\left(\frac{\pi^\e}{b}\right)+\frac{\pi^\e}{b}\nabla b
	\right]=H_b\left[\frac{\pi^\e}{b}\nabla b\right].
\end{align*}
    Thus, \eqref{Hbd} and \eqref{pi1} imply
\begin{align}\label{Pi}
	\|H_b[\n\pi^\e]\|_{L_t^\infty(H^2)}
	\leq C.
\end{align}
    The uniform bounds for $\n\psi^\e$ in $L^\infty(0,T;H^2(\Omega))\cap L^2(0,T;H^3(\Omega))$ and the Sobolev product estimates also give
\begin{align}\label{Ac}
	\|\dv(b\n\psi^\e\otimes\n\psi^\e)\|_{L_t^\infty(H^1)}+\|\dv(b\n\psi^\e\otimes\n\psi^\e)\|_{L_t^2(H^2)}
	\leq C.
\end{align}
    Consequently, we derive from \eqref{Hbd}, \eqref{Res}, \eqref{Pi}, and \eqref{Ac} that
\begin{align}\label{HR}
	\|H_b[\mathcal R^\e]\|
	_{L_t^\infty(H^1)}+\|\mathbb H_b[\mathcal R^\e]\|_{L_t^2(H^2)}\leq C.
\end{align}
    By \eqref{Hr} and \eqref{HR},
\begin{align}\label{lmt_5}
	H_b[\mathcal R^\e]\rightharpoonup 0\enspace
	{\rm in}\enspace L^2(0,T;H^2(\Omega)),\quad
	H_b[\mathcal R^\e]\overset{*}{\rightharpoonup}0 \enspace {\rm in}\enspace L^\infty(0,T;H^1(\Omega)).
\end{align}
    Combining \eqref{lmt_1}-\eqref{lmt_3}, \eqref{lmt_4}, and \eqref{lmt_5}, we pass to the limit in \eqref{Peq} and obtain
\begin{align}\label{Hlmt}
	H_b\left[(b\bv)_t+\dv(b\bv\otimes\bv)-\dv \mathbb{S}(\bv)\right]=0.
\end{align}
    By the weighted Helmholtz decomposition, there exists a scalar pressure $\Pi(t,x)$ such that
\begin{align*}
	(b\bv)_t+\dv(b\bv\otimes\bv)+b\n\Pi=\dv \mathbb{S}(\bv).
\end{align*}
    Together with $\dv(b\bv)=0$, $\bv\cdot\bn|_{\p\Omega}=0$, $\bn\times\crl\bv|_{\p\Omega}=0$ and $\bv|_{t=0}=\bv_0=b^{-1}H_{b}[b\mathbf{w}_0]$, this proves that $(\bv,\Pi)$ satisfies \eqref{LS_1}-\eqref{LS_2}. The proof of Theorem \ref{limit} is therefore finished. \qed
    
\subsection{Proof of Theorem \ref{conv}}
   The proof is motivated by the time-averaging technique, as used for example by Cheng \cite{Chen12} to establish convergence rates for the slow component of compressible Euler flows with solid-wall boundary condition and ill-prepared initial data. The key idea is to exploit the smallness of time primitives of the acoustic contributions and incorporate them into a corrected error. In this section, we combine \eqref{ac_1} and \eqref{lmt_4} with a weighted Neumann problem to rewrite the acoustic contributions to the projected equation as total time derivatives of $O(\e)$ corrections, up to $O(\e)$ remainders. Absorbing these corrections, including the acoustic contribution of the viscous term, into the projected momentum allows us to derive the desired convergence estimate by a weighted energy argument.
   
   Set $\mP^\e:=H_b[\rho^\e\bu^\e]$. By \eqref{rhop}, \eqref{unes}, \eqref{Hbd}, and Lemma \ref{PI},
\begin{align}\label{vel}
	&\|\phi^\e\|_{L_t^\infty(H^2)}+\|\e\phi^\e\|_{L_t^2(H^3)}+\|\bu^\e\|_{L_t^\infty(H^2)}+\|\bu^\e\|_{L_t^2(H^3)}+\|\pi^\e\|_{L_t^\infty(H^2)}+\|G^\e\|_{L_t^\infty(H^1)}\leq C,\\ \label{mom}
	&\|\mP^\e\|_{L_t^\infty(H^2)}+\|\mP^\e\|_{L_t^2(H^3)}+\|\n\psi^\e\|_{L_t^\infty(H^2)}+\|\n\psi^\e\|_{L_t^2(H^3)}\leq C.
\end{align}
    Projecting \eqref{eq1} gives
\begin{align*}
	\mathbf P_t^\e
	=H_b\left[\dv\mathbb S(\bu^\e)
	-\dv(\rho^\e\bu^\e\otimes\bu^\e)-\n\pi^\e\right],
\end{align*}
    and hence
\begin{equation}\label{Pt}
	\|\mathbf P_t^\e\|_{L_t^\infty(L^2)}
	+\|\mathbf P_t^\e\|_{L_t^2(H^1)}\leq C.
\end{equation}
    By \eqref{hydro} and conservation of mass, $\om\phi^\e dx=0$. We choose $\psi^\e$ defined in \eqref{neu} and define
\begin{align}\label{dchi}
	\chi^\e(t,x):=\frac{1}{\e}\int_0^t\psi^\e(s,x)ds.
\end{align}
    Integrating \eqref{ac_1} in time, we obtain
\begin{equation}\label{chi}
	\begin{cases}
		-\dv(b\n\chi^\e)=\phi^\e(t)-\phi^\e(0) & {\rm in}\enspace\Omega,\\
		\n\chi^\e\cdot\bn=0 & {\rm on }\enspace\p\Omega,\\
		\displaystyle\om\chi^\e\,dx=0.
	\end{cases}
\end{equation}
    Standard elliptic regularity for this  weighted Neumann problem gives
\begin{align}\label{nchi}
	\|\n\chi^\e\|_{L_t^\infty(H^3)}\leq C(\|\phi^\e\|_{L_t^\infty(H^2)}+\|\phi^\e(0)\|_{H^2})\leq C.
\end{align}
    The definition \eqref{dchi} also gives
\begin{equation}\label{primitive}
	\n\psi^\e=\e\p_t\n\chi^\e\quad{\rm in}\enspace L_t^\infty(H^2)\cap L_t^2(H^3).
\end{equation}
    Consequently,
\begin{align}\label{cross}
	H_b\left[\dv\left(\mP^\e\otimes\n\psi^\e+\n\psi^\e\otimes\mP^\e\right)\right]=&\e\p_tH_b\left[\dv\left(\mP^\e\otimes\n\chi^\e+\n\chi^\e\otimes\mP^\e\right)\right]\nonumber\\
	&-\e H_b\left[\dv\left(\mP_t^\e\otimes\n\chi^\e+\n\chi^\e\otimes\mP_t^\e\right)\right]\quad{\rm in}\enspace L_t^2(L^2).
\end{align}
    Since $\rho^\e\bu^\e=\mathbf P^\e+b\n\psi^\e$, we also have
\begin{align}\label{velocity}
	\bu^\e-\n\psi^\e
	=\frac{\mathbf P^\e}{b}
	-\frac{\e\phi^\e}{b}\bu^\e.
\end{align}
    The linearity of $\dv\mathbb S$ and \eqref{primitive} yield
\begin{align}\label{viscous}
	H_b[\dv\mathbb S(\bu^\e)]
	=&H_b\left[\dv\mathbb S\left(\frac{\mP^\e}{b}-\frac{\e\phi^\e}{b}\bu^\e\right)\right]+\e\p_tH_b[\dv\mathbb S(\n\chi^\e)].
\end{align}

    Define
\begin{align}\label{K}
	\mK^\e:=H_b\left[\dv(\mP^\e\otimes\n\chi^\e+\n\chi^\e\otimes\mP^\e)-\phi^\e\n\psi^\e-\dv\mathbb S(\n\chi^\e)\right].
\end{align}
    By \eqref{Hbd}, \eqref{vel}-\eqref{mom}, \eqref{nchi}, and Lemma \ref{PI}, the Sobolev product estimates then yield
\begin{align}\label{Kbound}
	\|\mK^\e\|_{L_t^\infty(H^1)}+\left\|\frac{\phi^\e}{b}\bu^\e\right\|_{L_t^\infty(H^2)}
	\leq C.
\end{align}
    Combining \eqref{Peq} with \eqref{cross}-\eqref{viscous}, we obtain
\begin{align}\label{equation}
	&\p_t(\mP^\e+\e\mK^\e)+H_b\left[\dv\left(\frac{\mP^\e\otimes\mP^\e}{b}\right)\right]\nonumber\\
	=&H_b\left[\dv\mathbb S\left(\frac{\mP^\e}{b}-\frac{\e\phi^\e}{b}\bu^\e\right)\right]-H_b\left[\dv\left(\left(\frac1{\rho^\e}-\frac1b\right)\rho^\e\bu^\e\otimes\rho^\e\bu^\e\right)\right]\nonumber\\
	&+\e H_b\left[\dv\left(\mP_t^\e\otimes\n\chi^\e+\n\chi^\e\otimes\mP_t^\e\right)\right]
	-\e H_b[\phi^\e F^\e+\n G^\e].
\end{align}
    The third terms on the right-hand side of \eqref{equation} can be estimated as follows:
\begin{align}\label{err3}
	&\e\left\|H_b\left[\dv\left(\mP_t^\e\otimes\n\chi^\e+\n\chi^\e\otimes\mP_t^\e\right)\right]\right\|_{L_t^2(L^2)}\leq C\|\mP_t^\e\|_{L_t^2(H^1)}\|\n\chi^\e\|_{L_t^\infty(H^3)}\leq C\e.
\end{align}
    where we have used \eqref{Pt} and \eqref{nchi}.
    
    To compare the corrected projected momentum $\mP^\e+\e\mK^\e$ with the limiting momentum $b\bv$, we introduce the corrected velocity error
\begin{align}\label{error}
	\me^\e:=\frac{\mP^\e+\e\mK^\e}{b}-\bv,
\end{align}
    we have
\begin{align*}
	\dv(b\me^\e)=0,\quad \me^\e\cdot\bn|_{\p\Omega}=0.
\end{align*}
    Thus, by the weighted orthogonality of $H_b$,
\begin{align*}
	\om H_b[f]\cdot\me^\e dx=\om f\cdot\me^\e dx.
\end{align*}
    Subtracting \eqref{Hlmt} from \eqref{equation} and using the definition of $\me^\e$, we obtain
\begin{align*}
	&b\p_t\me^\e
	+H_b\left[\dv\left(\frac{\mP^\e\otimes\mP^\e}{b}-b\bv\otimes\bv\right)\right]\\
	=&H_b\left[\dv\mathbb S(\bu^\e-\n\psi^\e-\bv)\right]-H_b\left[\dv\left(\left(\frac1{\rho^\e}-\frac1b\right)\rho^\e\bu^\e\otimes\rho^\e\bu^\e\right)\right]\\
	&+\e H_b\left[\dv\left(\mP_t^\e\otimes\n\chi^\e
	+\n\chi^\e\otimes\mP_t^\e\right)\right]-\e H_b[\phi^\e F^\e+\n G^\e].
\end{align*}
    
    For the viscous term, \eqref{velocity} gives
\begin{align*}
	\bu^\e-\n\psi^\e-\bv=\me^\e-\frac{\e\mK^\e}{b}-\frac{\e\phi^\e}{b}\bu^\e.
\end{align*}
    Since $\crl\n\psi^\e=0$ and $\n\psi^\e\cdot\bn|_{\p\Omega}=0$, we have $(\bu^\e-\n\psi^\e-\bv)\cdot\bn|_{\p\Omega}=0$ and $\crl(\bu^\e-\n\psi^\e-\bv)\times\bn|_{\p\Omega}=0$. By \eqref{Kbound} and Lemma \ref{dv-cl},
    integration by parts and Young's inequality yield
\begin{align}\label{errv}
	&\om\dv\mathbb S(\bu^\e-\n\psi^\e-\bv)
	\cdot\me^\e dx \nonumber\\
	\leq&-\mu\|\crl\me^\e\|_{L^2}^2-(2\mu+\lambda)\|\dv\me^\e\|_{L^2}^2+C\e\left\|\frac{\mathbf K^\e+\phi^\e\bu^\e}{b}\right\|_{H^1}
	\bigl(\|\crl\me^\e\|_{L^2}
	+\|\dv\me^\e\|_{L^2}\bigr) \nonumber\\
	\leq&-\frac{\mu}{2}\|\crl\me^\e\|_{L^2}^2
	-\frac{2\mu+\lambda}{2}\|\dv\me^\e\|_{L^2}^2
	+C\e^2 \nonumber\\
	\lesssim&-\|\me^\e\|_{H^1}^2
	+\|\me^\e\|_{L^2}^2+\e^2.
\end{align}
    For the convection terms, using $\dv\mP^\e=\dv(b\bv)=0$, we have
\begin{align}\label{pp}
	\dv\left(\frac{\mP^\e\otimes\mP^\e}{b}-b\bv\otimes\bv\right)
	=&(\mP^\e\cdot\n)\left(\frac{\mP^\e}{b}-\bv\right)+b\left[\left(\frac{\mP^\e}{b}-\bv\right)\cdot\n\right]\bv \nonumber\\
	=&(\mP^\e\cdot\n)\me^\e+b(\me^\e\cdot\n)\bv-\e\left[(\mP^\e\cdot\n)\frac{\mK^\e}{b}+b\left(\frac{\mK^\e}{b}\cdot\n\right)\bv\right].
\end{align}    
    Applying \eqref{mom}, \eqref{Kbound}, \eqref{pp}, and the boundary condition $\mP^\e\cdot\bn|_{\p\Omega}=0$, we obtain
\begin{align}\label{errc}
    &\left|\om\dv\left(\frac{\mP^\e\otimes\mP^\e}{b}-b\bv\otimes\bv\right)\cdot\me^\e dx\right| \nonumber\\
    \lesssim&\|\n\bv\|_{L^\infty}\|\sqrt{b}\me^\e\2+\e\|\mP^\e\|_{L^\infty}\left\|\n\left(\frac{\mK^\e}{b}\right)\right\|_{L^2}\|\me^\e\|_{L^2}+\e\|\mK^\e\|_{L^2}\|\n\bv\|_{L^\infty}\|\me^\e\|_{L^2} \nonumber\\
    \lesssim&\|\n\bv\|_{L^\infty}\|\sqrt{b}\me^\e\2+\e(1+\|\n\bv\|_{L^\infty})\|\me^\e\|_{L^2} \nonumber\\
    \lesssim&(1+\|\n\bv\|_{L^\infty})\|\sqrt{b}\me^\e\2+\e^2(1+\|\n\bv\|_{L^\infty}).
\end{align}

    Combining \eqref{err1}, \eqref{err2}, \eqref{err3}, \eqref{errv}, and \eqref{errc}, and taking $\delta$ sufficiently small, gives, for $0\leq t\leq T$,
\begin{align}\label{energy}
	&\|\sqrt{b}\me^\e(t)\2+\int_0^t\|\me^\e(s)\|_{H^1}^2ds\nonumber\\
	\leq&\|\sqrt{b}\me^\e(0)\2+C\e^2+C\int_0^t
	(1+\|\n\bv(s)\|_{L^\infty})\|\sqrt{b}\me^\e(s)\2ds.
\end{align}
    Since $\bv\in L_t^2(H^3)$,
\begin{align*}
	\int_0^T\|\n\bv(t)\|_{L^\infty}dt
	\leq CT^{1/2}\|\bv\|_{L_t^2(H^3)}<\infty.
\end{align*}
    Applying Gr\"onwall's inequality to \eqref{energy} and using \eqref{Kbound} and \eqref{error} yield
\begin{align*}
	\|\me^\e\|_{L_t^\infty(L^2)}^2
	+\|\me^\e\|_{L_t^2(H^1)}^2
	\leq C\left(\|\me^\e(0)\|_{L^2}^2+\e^2\right)\leq C\|H_b[\rho_0^\e\bu_0^\e]-b\bv_0\|_{L^2}+C\e.
\end{align*}
    Finally, using \eqref{error} and \eqref{Kbound} once more gives \eqref{unes}. \qed
	
	\vskip3mm
	\textbf{Acknowledgments} This work was partially supported by the National Natural Science Foundation of China under Grant Nos. 12131007 and 12471216. 
	
	\vskip2mm
	\textbf{Data availability}
	No datasets were generated or analyzed during the current study.

	\addcontentsline{toc}{section}{References}
	\bibliographystyle{plain} 
	\bibliography{ref} 

\end{document}